\documentclass[preprint,authoryear,sort&compress,11pt]{elsarticle}
\pdfoutput=1
\usepackage{fullpage}
\usepackage[utf8]{inputenc}
\usepackage[english]{babel}
\usepackage{graphicx,amsmath,amssymb, amsthm, multicol, bm}
\usepackage{url} 
\usepackage[hidelinks,unicode,verbose,hypertexnames=false]{hyperref}
\usepackage[dvipsnames]{xcolor}
\usepackage{tablefootnote}
\usepackage[shortlabels]{enumitem}
\usepackage{multirow}
\usepackage{array}
\usepackage{txfonts}
\usepackage{tocbibind}
\usepackage{pifont}
\usepackage{float}
\renewcommand{\checkmark}{\ding{51}}
\newcommand{\xmark}{\ding{53}}

\makeatletter
\newcommand{\myitem}[1]{%
\item[#1]\protected@edef\@currentlabel{#1}%
}
\makeatother

\journal{European Journal of Operational Research}
\makeatletter
\def\ps@pprintTitle{%
     \let\@oddhead\@empty
		 \let\@evenhead\@empty
     \def\@oddfoot
       {\hbox to \textwidth%
        {\ifnopreprintline\relax\else
        \@myfooterfont%
         \ifx\@elsarticlemyfooteralign\@elsarticlemyfooteraligncenter%
           \hfil\@elsarticlemyfooter\hfil%
         \else%
         \ifx\@elsarticlemyfooteralign\@elsarticlemyfooteralignleft%
           \@elsarticlemyfooter\hfill{}%
         \else%
         \ifx\@elsarticlemyfooteralign\@elsarticlemyfooteralignright%
           {}\hfill\@elsarticlemyfooter%
         \else%
            \begin{minipage}[t]{\textwidth}   Preprint accepted in \ifx\@journal\@empty%
                 Elsevier%
            \else\@journal\fi\hfill June 28, 2023\\
						Published version available at \rurl{10.1016/j.ejor.2023.06.039}\medskip\\ 
						\begin{minipage}[b]{0.7\textwidth}\textcopyright\ 2023. This work is licensed under a Creative Commons\\
			      Attribution-NonCommercial-NoDerivatives
						4.0 International License\end{minipage} 
							\hfill\includegraphics[width=0.2\textwidth]{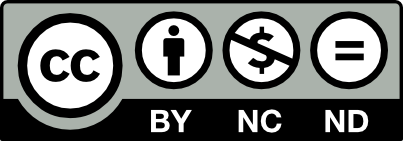}
						\end{minipage}\fi%
         \fi%
         \fi%
         \fi%
         }
       }%
     \let\@evenfoot\@oddfoot}
\makeatother
\newcommand\rurl[1]{%
  \href{http://doi.org/#1}{\nolinkurl{#1}}%
}

\newtheorem{theorem} {Theorem}

\newtheorem{lemma} {Lemma}

\newtheorem{definition} {Definition}

\theoremstyle{remark}
\newtheorem{pozn}{Remark}

\newcommand{\DMU}{\mathrm{DMU}}
\newcommand{\ones}{\bm{e}}
\newcommand{\T}{\mathcal{T}}
\newcommand{\dom}{\operatorname{dom}}
\newcommand{\im}{\operatorname{im}}
\newcommand{\R}{\mathbb{R}}
\newcommand{\m}{\mathrm{min}}
\newcommand{\M}{\mathrm{max}}
\newcommand{\X}{\bm{X}}
\newcommand{\Y}{\bm{Y}}
\newcommand{\x}{\bm{x}}
\newcommand{\y}{\bm{y}}
\newcommand{\g}{\bm{g}}
\newcommand{\p}{\bm{\phi}}
\newcommand{\la}{\bm{\lambda}}
\newcommand{\0}{\bm{0}}
\newcommand{\s}{\bm{s}}
\newcommand{\e}{\mathrm{e}}
\newcommand{\D}{\mathcal{D}}
\makeatletter
\newcommand\newtag[2]{#1\def\@currentlabel{#1}\label{#2}}
\makeatother

\renewenvironment{thebibliography}[1]{%
\begin{oldthebibliography}{#1}%
\setlength{\baselineskip}{1.2em}
\linespread{1.2}
\setlength{\parskip}{.3ex}%
\setlength{\itemsep}{.3em}%
}%
{%
\end{oldthebibliography}%
}

\begin{document}

\begin{frontmatter}

\title{
A unified approach to radial, hyperbolic, and directional efficiency measurement in Data Envelopment Analysis
}
\author[fmfi]{Margar\'eta Halick\'a}
\ead{halicka@fmph.uniba.sk}
\author[fmfi]{M\'aria Trnovsk\'a\texorpdfstring{\corref{cor1}}{}}
\ead{trnovska@fmph.uniba.sk}
\author[bayes]{Ale\v{s} \v{C}ern\'{y}}
\ead{ales.cerny.1@city.ac.uk}

\cortext[cor1]{Corresponding author}
\address[fmfi]{Faculty of Mathematics,  Physics and Informatics, Comenius
 University in Bratislava,  Mlynsk\'a dolina, 842 48 Bratislava, Slovakia}
 \address[bayes]{Bayes Business School, City, University of London, 106 Bunhill Row, London EC1Y 8TZ, UK}
\date{April 20, 2023}
\begin{abstract} 
The paper analyses properties of a large class of ``path-based'' Data Envelopment Analysis models through a unifying general  scheme. The scheme includes the well-known oriented radial models, the hyperbolic distance function model,  the directional distance function models, and even permits their generalisations.  The modelling is not constrained to non-negative data and is flexible enough to accommodate variants of standard models over \emph{arbitrary data}.

Mathematical tools developed in the paper allow systematic analysis of the models from the point of view of ten desirable properties.  It is shown that some of the properties are satisfied (resp., fail) for all models in the general scheme, while others have a more nuanced behaviour and must be assessed individually in each model. Our results can help researchers and practitioners navigate among the different models and apply the models to mixed data.
\end{abstract}

\end{frontmatter}
\setlength{\abovedisplayskip}{4pt}
\setlength{\belowdisplayskip}{4pt}

\section{Introduction}

Data Envelopment Analysis (DEA) is a non-parametric analytical method used to assess the performance of Decision Making Units (DMU).
The DEA approach defines technology sets via observed input-output data in combination with certain axioms.  Each traditional model links an efficiency measure and a given technology set using a mathematical optimisation programme.  By solving the programme for an assessed DMU, one finds both the efficiency score and a projection%
\footnote{Synonymously with ``projection,'' the literature also uses the terminology ``reference unit,'' ``projection benchmark,'' ``projection point,'' or ``target'' in this context.}
on the frontier of the technology set that dominates the assessed DMU.%
\footnote{The resulting projection need not be Pareto--Coopmans efficient.}
\citet{russell.schworm.18} distinguish between two main methods of arriving at the efficiency score, giving rise to two classes of DEA models appositely named \emph{slacks-based} and \emph{path-based} models. 

This paper establishes a \emph{general path-based scheme} built on variable returns to scale (VRS) technology that provides a unifying mathematical frame\-work for the well-known oriented radial models, the hyperbolic distance function model, and the directional distance function models. The scheme offers a rich menu of projection paths for efficiency analysis and the special capability of handling negative data, which is important in a variety of applications, especially in the areas of accounting and finance.  In the confines of the general scheme, we provide a systematic and comprehensive analysis of the models \emph{vis-\`a-vis} ten desirable properties. Overall, our paper amplifies the message of \citet{russell.schworm.18} that recognising the class to which a model belongs allows one to deduce some of the model characteristics, which is important for the correct interpretation of models and their practical use.

The slacks-based and the path-based models are easily distinguishable from each other on the basis of the objective function in the corresponding programmes, or, in other words, based on the way they look for a projection.
According to \citet{russell.schworm.18}, the models in the \emph{slacks-based class}  search for the projection
by ``specifying the form of aggregation over the coordinate-wise slacks.'' The slacks indicate the input surplus and output shortage between the projection and the assessed DMU. 
The main representatives of this class are the Slacks-Based Measure (SBM) model of \cite{tone.01}, the Russell Graph Measure model of \cite*{fare.al.85}, the Additive Model (AM) of \cite{charnes.al.85}, and the Weighted Additive Models (WAM) including the Range Adjusted Measure (RAM) model  (\citealp*{cooper.al.99}) and  the Bounded Adjusted Measure (BAM) model (\citealp{cooper.11}).

The \emph{path-based models} --- the main focus of this paper --- search for the projection by specifying various parametric paths which run from the assessed DMU to the boundary of the technology set.  In the special case of the radial BCC input or output-oriented models (\citealp*{banker.al.84}), the path is defined by a ray connecting the DMU to the origin in the space of inputs or outputs and thus represents a proportional, radial contraction of inputs or expansion of outputs. In the Directional Distance Function (DDF) model (\citealp*{chambers.al.96.jet,chambers.al.98}), the path is determined by a ray in a pre-assigned  direction pointing from the assessed DMU towards  the dominating part of the frontier.%
\footnote{There are also approaches, where the directional vector may not point to the dominating part of the frontier, or may even be determined endogenously. For a discussion on these non-standard approaches we refer to \cite{pastor.al.22}.}
On this path, one then seeks the point of minimal distance to the frontier of the technology set. The Hyperbolic Distance  Function (HDF) model (\citealp*{fare.al.85}) combines together  the  input and output oriented BCC models by using a hyperbolic path that allows for a simultaneous equiproportionate contraction of inputs and expansion of outputs.

Despite the many papers on DEA, there are only a few studies analysing the properties of either class of models 
in a unified framework. In the context of general economic productivity theory, a series of articles \cite*{russell.schworm.08}, \cite*{russell.schworm.11}, \cite*{levkoff.al.11}, \cite*{roshdi.al.18}, and \cite*{hasannasab.al.19} provided a comprehensive analysis of efficiency measures  over different types of productivity sets. The main message  from these papers is  that the slacks-based measures, when operating in the full input-output space, identify the Pareto--Koopmans efficiency unambiguously while the path-based measures do not. 

In the DEA setting,
 \cite*{cooper.al.99} and \cite{sueyoshi.sekitani.09.ejor.196} analyse efficiency measures and introduce a set of desirable properties that an ideal DEA model should satisfy. These properties include the indication of strong efficiency; boundedness; strict monotonicity;  unit invariance; and translation invariance. The selected DEA models are then classified on the basis of these criteria.

In other work, \cite{halicka.trnovska.21} analyse the properties of slacks-based models in a general scheme that encompasses all commonly used models in this class, but also allows for the construction of new models.%
\footnote{The terminology ``non-radial'' and ``radial pattern models'' was used for slacks-based and path-based models, respectively.}
  All models in the general slacks-based scheme project
  onto the strongly efficient frontier and therefore account for all sources of inefficiencies. Among them, the RAM model performs best when measured against eight desirable properties, satisfying seven, and failing only the unique projection property due to its linearity. Recently, \citet{aparicio.monge.22} have proposed a convex generalisation of the RAM model that falls into the general scheme of \cite{halicka.trnovska.21} and provides a unique projection point.  Thus, the new model currently claims the top spot in terms of the number of desirable properties.

Our paper aims to redress the lack of comprehensive analysis of path-based models in the literature.  The knowledge about the path-based models is currently fragmented across many articles with varying focus. This situation is exacerbated by the fact that the properties of, for example, DDF models depend significantly on the choice of direction vectors and these have not been treated systematically to date. 
There appears to be a general consensus that the path-based models do not guarantee strongly efficient projection points and, therefore, their efficiency score is overstated; and that they are monotone but not strictly monotone. 
Many authors noticed difficulties with super-efficiency measurement under variable returns to scale and the associated measurement of productivity change over time (e.g., \citealp{briec.kerstens.09}).
\cite{aparicio.al.16} investigated the translation invariance of DDF models.
A certain type of homogeneity was observed in the oriented radial models and the HDF model (\citealp{cuesta.zofio.05}). 
 On the other hand, DDF models have the property of homogeneity only in the case of constant returns to scale.  Another type of homogeneity (so-called g-homogene\-ity) was introduced to describe some properties of DDF  (\citealp{hudgins.primont.07}).

 In this paper, we analyse path-based models in a general framework. The main building block of our approach is a parametric path that starts at the assessed DMU and for decreasing values of the parameter passes through dominating units in the technology set towards its frontier. The path has two main ingredients: (i) a direction vector; (ii) a real-valued smooth function, whose specifications determine the path shape and, ultimately, the model properties.  The assumptions imposed on directions and path functions are flexible enough to accommodate all standard models, such as the BCC input and output-oriented models, the HDF and DDF models, the generalised distance function model by \cite{chavas.cox.99}, and even offer the possibility of going beyond.  
 
 A further advantage of the proposed general framework is that it permits extension of existing models designed for non-negative  data to arbitrary data (i.e., it accommodates the input and output data, for which some or all components are negative). 
 As a rule, DEA models are designed for non-negative data. In practice, however, 
negative data
 are encountered in many applications in areas such as insurance, accounting, finance, or banking. A common approach to overcoming this difficulty is to use translation-invariant models, which can be applied to arbitrary data without modifications.  The drawback is that only a few standard models have the property of translation invariance.  
 These include some of the DDF models (see \citealp{aparicio.al.16}) but not BCC or HDF.  Therefore, many DEA studies propose procedures of varying complexity to deal with negative data (e.g., \citealp{cheng.al.13}; \citealp{tone.al.20}). 

The general framework for path-based models allows us to survey in one place the properties of all standard path-based models and also offer certain guidance on how to construct new models with given properties. The desirable properties analysed in this paper include (a) unique projection point; (b) indication; (c) strong efficiency of projection points; (d) boundedness; (e) unit and translation invariance; (f) (strict) monotonicity; (g) super-efficiency;  and (h) homogeneity.  Our results divide the properties into two groups: the properties that hold universally (resp., universally fail)  
for every model in the general framework; and the properties that must be assessed individually for each model. With the help of the general framework, we show that only the unique projection point property is satisfied in all models.  Three other properties (indication, strong efficiency of projection, and strict monotonicity) usually fail simultaneously, although surprisingly there are very special cases where they simultaneously hold.  For properties in the second group, the article provides tools to determine whether a property holds or fails in a specific model. In particular, the monotonicity property is satisfied by all standard path-based models, but the homogeneity property is limited to only a specific type of directions and path functions. We find that there are no trade-offs between the properties of indication, strong efficiency, and strict monotonicity on the one hand and the property of homogeneity on the other.

The paper is organised as follows. Section \ref{S2} introduces basic terminology and notation concerning, among others, the technological set and its efficient frontier over general data; desirable properties of DEA models;   and standard path-based models. Section \ref{S3} proposes a general scheme for path-based models, analyses its basic properties, and discusses its geometry. Section \ref{S4} conducts a deeper analysis of the general scheme in light of ten desirable properties. This section develops practical criteria for each property and illustrates them on individual standard path-based models.
Section \ref{S5} extends the above analysis in two directions: (a) properties of standard models over arbitrary data;
(b) construction of new models with good properties. Section \ref{S6} concludes.

\section{Preliminaries}\label{S2}
Let us establish some notation. $\R^d$ denotes the $d$ dimensional Euclidean space, and $\R^d_+$ its non-negative orthant. The lowercase bold letters denote column vectors 
and the uppercase bold letters matrices.
The superscript $T$ denotes the transpose of a column vector or a matrix.
The symbol $\ones$ denotes a vector of ones.

Consider a set of $n$ observed decision-making units ${\DMU}_j\ (j=1,\dots, n)$, each consuming $m$ inputs $x_{ij}$ $(i=1,\dots,m)$ to produce $s$ outputs $y_{rj}\ (r=1,\dots,s)$. For each $j=1,\dots,n$, the data of the inputs and outputs of ${\DMU}_j$ can be arranged into the column vectors $\x_j=(x_{1j},\dots,x_{mj})^T\in \R^m$ of the inputs and $\y_j=(y_{1j},\dots,y_{sj})^T\in \R^s$ of the outputs.  Finally, the input and output vectors of all DMUs form the $m\times n$ input and $s\times n$ output matrices $\X$ and $\Y$, i.e.,  $\X=[\x_1,\dots,
\x_n]$  and $\Y=[\y_1,\dots ,\y_n]$,  respectively.

In the article, we reserve the notations $i$, $j$, and $r$ for the indices that go through whole index sets $\{1,\dots,m\}$, $\{1,\dots,n\}$, and $\{1,\dots,s\}$, respectively. We make no assumptions about the non-negativity of the data at this point. The non-negativity requirement may follow later from other assumptions, and we shall alert the reader whenever that is the case.

\subsection{Technology set}\label{SS:T}

On the basis of the given data, we consider the following technology set:
\begin{equation}\label{T}
{\T}=\left\{(\x,\y)\in\R^m\times \R^s\ | \
\X\la\le \x,\ \Y\la\ge \y,\ \la\ge \0, \ \ones^T\la =1
\right\},
\end{equation}
corresponding to variable returns to scale (VRS). Note that the common non-negativity of $(\x,\y)$ is not imposed here. Elements of ${\T}$ will be called units. It follows from \eqref{T} that the closed set $\T$ has a non-empty interior; we shall denote its boundary by $\partial \T:=\T\setminus \text{int}\T$.

By $(\x_o,\y_o)$ we denote a unit from ${\T}$ to be currently
evaluated. For input vectors, we also use the notation $\x^\m,\x^\M,\x^{ev},\x^{sd}$,
where for $i=1,\dots m$ we set $x^{\m}_i=\min_jx_{ij}$, $x^{\M}_i=\max_jx_{ij}$, $x^{ev}_i=\frac 1n\sum_j x_{ij}$,  and
$x^{sd}_i = \sqrt{\frac1n\sum_j (x_{ij}-\x^{ev}_i)^2}.$
Here $x^{sd}_i$ is the standard deviation of the $i$-th input over all $\DMU_j$, $j=1,\dots,n$.
The notation $\y^\m$, $\y^\M$, $\y^{ev}$, and $\y^{sd}$  is introduced analogously for the output vectors.  Without loss of generality, we assume that $\x^\m<\x^\M$ and $\y^\m<\y^\M$. Otherwise, there would be components of inputs / outputs, where all DMUs take the same value, and such components can be excluded from the analysis.

We write $(\x,\y)\succsim (\x_o,\y_o)$ if the unit $(\x, \y)$ \emph{dominates} the unit $( \x_o, \y_o)$, that is, if $ \x\le \x_o$ and $ \y\ge\y_o$. 
  A unit $( \x, \y)$ \emph{strictly dominates} the unit $( \x_o, \y_o)$ if $ \x<\x_o$ and $ \y> \y_o$.
  A unit $( \x_o, \y_o)\in \T$ is called \emph{strongly efficient} if there is no other unit in $\T$ that dominates $( \x_o, \y_o)$, that is, if $( \x, \y)\in \T$  dominates $( \x_o, \y_o)$, then $( \x, \y)=( \x_o, \y_o)$.%
  \footnote{This is the well known Pareto--Koopmans efficiency.   Some authors call such units Pareto efficient, or fully efficient; see the discussion in \citet[p. 45]{cooper.al.07}.}
  A unit $(\x_o,\y_o)\in \T$ is called \emph{weakly efficient} if there is no unit in $\T$ that strictly dominates $(\x_o,\y_o)$.
Evidently, any strongly efficient unit is weakly efficient, and weakly efficient units lie on the boundary $\partial\T$.

The converse is also true: every unit on the boundary $\partial \T$ is weakly efficient because the definition of $\T$ in \eqref{T} does not impose the non-negativity assumption on the units therein. Therefore, the boundary $\partial \T$ is partitioned into the \emph{strongly efficient frontier} $\partial^S\T$ containing all strongly efficient units and the remaining part $\partial^W\T:= \partial\T\setminus \partial^S\T$, which consists of weakly but not strongly efficient units. In this paper, we refer to the remaining part of the boundary as \emph{weakly efficient frontier}. One thus has $\partial \T=\partial^S\T\cup\partial^W\T$ and $\partial^S\T\cap\partial^W\T=\emptyset$. The simple proof of the next lemma is omitted.

\begin{lemma}\label{L:weak} For $(\x_o,\y_o)\in\T$, the following statements hold.
\begin{enumerate}[(a)]
 \item\label{W.o}
 $(\x^{\m},\y^{\M})\succsim(\x_o,\y_o)$;
 \item\label{W.a}
 $(\x_o,\y_o)\in\partial\T$ if and only if for all $(\bm d^x,\bm d^y)\ge\0$
  such that $( \x_{o}-\bm d^x, \y_o+\bm d^y)\in \T$, one has $(\bm d^x,\bm d^y)\ngtr\0$;
 \item\label{W.b}
  $(\x_o,\y_o)\in\partial^S\T$ if and only if for all $(\bm d^x,\bm d^y)\ge\0$ such that $( \x_{o}-\bm d^x, \y_o+\bm d^y)\in \T$, one has $(\bm d^x,\bm d^y)=\0$;
 \item\label{W.c}
  $(\x_o,\y_o)\in\partial^W\T$ if and only if there exists $(\bm d^x,\bm d^y)\gneqq\0$ such that $( \x_{o}-\bm d^x, \y_o+\bm d^y)\in \T$ and if for all $(\bm d^x,\bm d^y)\ge\0$
  such that $( \x_{o}-\bm d^x, \y_o+\bm d^y)\in \T$, one has $(\bm d^x,\bm d^y)\ngtr\0$.
  \end{enumerate}
\end{lemma}

\subsection{Efficiency measures and their desirable properties}

On the technology set $\T$, one can define various \emph{measures}, which then give rise to a specific DEA \emph{models}. This is achieved by formulating certain mathematical programming problems that are applied separately to each evaluated unit $(\x_o,\y_o)\in\T$.

A mathematical programme consists of an \emph{objective function} and a set of conditions for \emph{decision variables}. The values of the decision variables that satisfy these conditions are called \emph{feasible solutions} of the programme. Solving the programme yields the \emph{optimal value} of the objective function and a subset of feasible solutions called \emph{optimal solutions}. 
The optimal value of a DEA model determines the \emph{efficiency score / value of the measure} of the evaluated unit.  Optimal solutions also yield virtual units  called \emph{projections}  of the evaluated unit ---  the elements of $\partial\T$ that dominate $(\x_o,\y_o)$ and are associated with the efficiency score.

Depending on the chosen measure, models and their efficiency scores have varying properties. Following the work of \cite{fare.lovell.78}, \cite{pastor.al.99}, \cite{sueyoshi.sekitani.09.ejor.196}, \citet{russell.schworm.18}, \cite{aparicio.monge.22}, and
\cite{halicka.trnovska.21}, we highlight ten desirable properties for efficiency measurement in DEA.
\begin{enumerate}[label={\rm(P\arabic{*})}, ref={\rm(P\arabic{*})}]
\item\label{P1}  {\bf Unique projection for efficiency comparison.} While the model may admit multiple 
optimal solutions, the projection point is unique.

\item\label{P2} {\bf Indication.} The value of the measure equals one if and only if the evaluated unit is strongly efficient.
\item\label{P3} {\bf Strong efficiency of projections.}  Projections generated by the measure are strongly efficient.
\item\label{P4} {\bf Boundedness.} The measure takes values between zero and one.
\item\label{P5}  {\bf Units invariance.} The value of the measure does not depend on the units of measurement in the input and output variables.
\item\label{P6} {\bf Translation invariance.} The value of the measure
is not affected by translation of inputs or outputs.
\item\label{P7}  {\bf Monotonicity.} An increase in any input or a decrease in any output relative to the evaluated unit, keeping other inputs and outputs constant, reduces or maintains the value of the measure.
\item\label{P8}  {\bf Strict monotonicity.} An increase in any input or a decrease in any output relative to the evaluated unit, keeping other inputs and outputs constant, reduces the value of the measure.
\item\label{P9} {\bf Super-efficiency.} The value of the measure of a unit outside the technology set is well defined and finite.
\item\label{P10} {\bf Homogeneity.} Feasible scaling of the input vector and the output vector of the evaluated unit by a power of $\mu>0$ (the power may be different for inputs and outputs) results in the efficiency score being scaled by a power of $\mu$.
\end{enumerate}

The importance of individual properties varies with  the goals of the analysis. In the first instance, one needs to know whether the efficiency score equal to one indicates the strong efficiency of the evaluated unit. Therefore, \ref{P2} is perhaps the most important property. It is closely related to \ref{P3}, although the two properties are not equivalent. Failure of \ref{P3} indicates that the score may be overestimated.   
Some models are ambiguous  in terms of satisfying the properties \ref{P2} and \ref{P3} and require further analysis to reach a conclusion.

Invariance properties \ref{P5} and \ref{P6} are important when the input and output data are expressed in different measurement units or when the data also contain negative values, respectively. 
No universal conclusions can be drawn regarding these properties, hence each DEA model must be assessed individually. 
The property \ref{P6} is tied to variable returns-to-scale technology set.

The monotonicity properties \ref{P7} and \ref{P8} ensure that competing units are evaluated with  a certain degree of fairness, in that a dominating unit receives a score not below the score of a dominated unit. 
The weak monotonicity \ref{P7} is usually  satisfied in standard DEA models and the strong monotonicity \ref{P8} is satisfied in most slacks-based models.

The property \ref{P9} is important for (i) measuring super-efficiency, which allows one to discriminate among efficient units; and (ii) measuring the frontier shift between
different time periods; see \cite{sueyoshi.sekitani.09.ejor.196}. 
Considerable difficulties with this property  relate particularly to VRS, where  some super-efficiency models are generally not feasible.

The homogeneity property  \ref{P10}, together with  \ref{P2},  \ref{P3}, and  \ref{P8}, was formulated for the first time as a desirable property of efficiency measures in \cite{fare.lovell.78} in connection with input-oriented measures. Since then, it has been subject to several generalisations suitable for output-oriented models and graph models, which is also reflected in the terminology  (e.g., ``almost homogeneity''). For slacks-based models, the concept was relaxed to sub-homogeneity (\citealp{pastor.al.99}), and for the directional distance function, it was modified to g-homogeneity (\citealp{hudgins.primont.07}).

The property of boundedness \ref{P4} was included in the list of desirable properties in \cite {fare.al.83}.  The property \ref{P4} guarantees that the efficiency measure is bounded by 0
and 1, where the value 1 identifies the strong efficiency in \ref{P2}.

In contrast to all the remaining properties that are expressed in terms of the efficiency measure, the properties \ref{P1} and  \ref{P3} are expressed using the projections provided by the models. The importance of \ref{P1} was highlighted in \cite{sueyoshi.sekitani.09.ejor.196}, especially in connection with the slacks-based models, where the concept of projection is not uniquely defined.

It is well known that no model meets all properties \ref{P1}--\ref{P10}. However, some of the properties occur in groups and hold simultaneously for a class of models.
\cite{halicka.trnovska.21} show that slacks-based models satisfy \ref{P2}, \ref{P3}, \ref{P7}, \ref{P8},   and violate \ref{P1}. The present paper shows that path-based models satisfy \ref{P1} and \ref{P7}, and fail \ref{P2}, \ref{P3}, and \ref{P8}.

We shall now introduce some standard models that enter our analysis as special cases. We formulate these models in a way that allows for certain generalisations.

\subsection{Standard path-based models} First, we present  input/output-oriented radial models, which are known in DEA under the BCC acronym  in honour of their authors \cite{banker.al.84}, who introduced them to DEA in the VRS version. Note that the standard input/output-oriented DEA radial models are inspired by the Shephard’s  (\citealp{shephard.53}) 
input/output distance functions, and the technical efficiency measurement is based on the seminal work of \cite{farrell.57}. The standard DEA formulation of these models is over non-negative data with the assumption
\begin{equation}\label{xygneqq0}
\x_j\gneqq \0 \quad \text{ and }\quad  \y_j\gneqq \0\qquad \text{ for each } j\in\{1,\dots,n\}.
\end{equation}
The input oriented BCC model can be formulated as
\begin{equation}\label{BCCI}
\begin{array}{rl}
\text{BCC-I:}\qquad\min & \theta \\
     & \X\la \le \theta \x_o \\
     & \Y\la \ge \y_o \\
     & \ones^T\la =1,
     \ \la \ge \0,
\end{array}
\end{equation}
and the output oriented BCC model
as
\begin{equation}\label{BCCO}
\begin{array}{rl}
\text{BCC-O:}\qquad\min & \theta \\
     & \X\la \le  \x_o \\
     & \Y\la \ge \frac{1}{\theta} \y_o \\
     & \ones^T\la =1, \
      \la \ge \0.
\end{array}
\end{equation}
If we denote $\psi:=\frac{1}{\theta}$ and replace $\min \theta$ by $\max \psi$, then we obtain the output BCC model in standard form, where the optimal $\psi^*\ge 1$ and the efficiency score is $\frac{1}{\psi^*}$.
Among the properties \ref{P1}--\ref{P9}, the BCC models satisfy only \ref{P1}, \ref{P4}, \ref{P5}, and \ref{P7}.

The hyperbolic distance function model introduced by \cite{fare.al.85} combines the input and output oriented radial measures into one measure.
Commonly, it is defined over positive data only, i.e.
\begin{equation}\label{xy>0}
\x_j> \0 \quad \text{ and }\quad  \y_j> \0\qquad \text{ for each } j\in\{1,\dots,n\}.
\end{equation}
In our formulation, it reads:
\begin{equation}\label{HDF}
\begin{array}{rl}
\text{HDF:}\qquad\min& \theta \\
      &\X\la \le \theta \x_o \\
      &\Y\la \ge \frac{1}{\theta} \y_o \\
      &\ones^T\la=1, \
      \la \ge \0.
\end{array}
\end{equation}
\cite{halicka.trnovska.19} investigated the properties of this model, proposed computational methods for its solution, and derived its dual form through SDP, while \cite{hasannasab.al.19} established duality through the second-order cone. HDF satisfies the properties \ref{P1}, \ref{P4}, \ref{P5},  \ref{P7}, and \ref{P9}.

 The use of linearisation $\frac{1}{\theta}-1\sim 1-\theta $ around $\theta=1$ in the HDF model (see \citealp{fare.al.16}) leads to the following:
\begin{equation}\label{DDF}
\begin{array}{rl}
\text{DDF:}\qquad\min & \theta \\
     & \X\la \le \x_o- (1-\theta) \x_o \\
     & \Y\la \ge \y_o +(1-\theta)\y_o\\
     & \ones^T\la =1, \
      \la \ge \0.
\end{array}
\end{equation}
If we denote $\delta=1-\theta$ and replace $\min \theta$ by $\max \delta$, we obtain a special version of the directional distance function model with the directional vector $(g_o^x, g_o^y)=(x_o,y_o)$.  This model, introduced by \cite{briec.97}, is known as the proportional directional distance function. It satisfies the properties \ref{P1}, \ref{P4}, \ref{P5}, \ref{P7}, and \ref{P9}.

Now, we present the general directional distance function model (DDF-g). Let $ \g_o^x\in \R^m_+$, $\g_o^y\in \R^s_+$ be given directional vectors that may depend on $(\x_o,\y_o)$ and such that at least one component of $ \g_o^x$ or $\g_o^y$ is positive, that is,
\begin{equation}\label{g>neq0}
(\g_o^x,\g_o^y)\gneqq \0.
\end{equation}
The model admits negative data, provided that the assumption of non-negativity \eqref{g>neq0} on the directional vectors is satisfied. The model reads as follows:
\begin{equation}\label{DDF-G}
\begin{array}{rl}
\text{DDF-g:}\qquad\min & \theta \\
     & \X\la \le \x_o- (1-\theta) \g_o^x\\
     & \Y\la \ge \y_o +(1-\theta)\g_o^y\\
     & \ones^T\la =1, \
      \la \ge \0.
\end{array}
\end{equation}
Once again, if we denote $\delta=1-\theta$ and replace $\min \theta$ by $\max \delta$, we obtain the standard form of DDF-g  introduced in \cite{chambers.al.96.jet,chambers.al.98}.

Since DDF is a linearisation of HDF, it is natural to ask whether there exists a generalisation of HDF whose linearisation corresponds to DDF-g. An affirmative answer is provided by the following model, which we refer to as the general hyperbolic distance function (HDF-g) model:
\begin{equation}\label{HDF-G}
\begin{array}{rl}
\text{HDF-g:}\qquad\min & \theta \\
     & \X\la \le \x_o- (1-\theta) \g_o^x\\
     & \Y\la \ge \y_o +(\frac{1}{\theta}-1)\g_o^y\\
     & \ones^T\la =1, \
      \la \ge \0.
\end{array}
\end{equation}
This new model, too, uses  pre-specified directions and admits negative data as long as the directions satisfy the non-negativity condition \eqref{g>neq0}. It is easy to see that the choices $\g_o^x=\x_o$ and $\g_o^y=\y_o$ lead to the HDF model.

Without specific knowledge of the directional vector $\g$ beyond \eqref{g>neq0}, only the property \ref{P1} is ensured in HDF-g and DDF-g. More must be said about the choice of directional vectors to ensure the validity of the remaining properties \ref{P2}--\ref{P9}.
For a discussion of the available choices of directional vectors, see, for example, \cite{fare.al.08} and \cite{pastor.al.22}.
Table~\ref{Tab:choices} presents the directions $(\g_o^x,\g_o^y)$ that will be included in the subsequent analysis of the families of HDF-g and DDF-g models. These directions exhibit a range of properties, as documented in Table~\ref{Tab:properties}.
\begin{table}[t!]
\centering
 \begin{tabular}{l c c c}
 \hline\\[-1.8ex]
 Notation  & $\g_o^x$  & $\g_o^y$ & Reference\\[0.8ex]
 \hline\\[-1.8ex]
 \newtag{(G1)}{G1} & $\x_o$  & $\y_o$ & \cite{chambers.al.96.per} \\
 \newtag{(G2)}{G2} & $\x_o- \x^{\m}$ & $ \y^{\M}-\y_o$ &\cite{portela.al.04} \\
 \newtag{(G3)}{G3}&$ \x^{\M}- \x^{\m}$&$ \y^{\M}- \y^{\m}$&\cite{portela.al.04}\\
 \newtag{(G4)}{G4}&$\x^{ev}$&$\y^{ev}$&\cite{aparicio.al.13}\\
 \newtag{(G5)}{G5}&$\x^{sd}$&$\y^{sd}$&\\
 \newtag{(G6)}{G6}&$\ones=(1,\dots,1)^T$&$\ones=(1,\dots,1)^T$& \cite{chambers.al.96.per}\\[0.8ex]
 \hline
 \end{tabular}
 \caption{Commonly encountered choices of directions $g_o$. For the notation $\x^{\m}$, $\x^{\M}$, etc., see Subsection~\ref{SS:T}.}
\label{Tab:choices}
\end{table}

\ref{G1} is perhaps the most widely chosen direction in the family of DDF-g models, which yields the Farrell proportional distance function (see \citealp{briec.kerstens.09}). The directions \ref{G2}  and \ref{G3} appear in the range-directional models developed by \cite{portela.al.04}. The use of \ref{G4} directions is discussed in \cite{aparicio.al.13}.
 It is also quite common to consider DDF-g with \ref{G6} directions, which is mathematically equivalent to minimising the $l_{\infty}$ distance to the boundary of $\T$ (see \citealp{briec.lesourd.99}).

\begin{pozn}It is interesting to note that the expressions for directions \ref{G1}--\ref{G6} appear in the denominator of weights in specific WAM  models. Indeed, direction \ref{G6} is linked to the weights in AD model; directions  \ref{G2}  and \ref{G3} are linked to the weights in BAM and RAM models, respectively; direction \ref{G1} corresponds to the weights in so-called Measure of Ineficiency Proportions (MIP) model by \cite{cooper.al.99}; and the weights generated via \ref{G5} appear in the WAM model of \cite{sevcovic.al.01}.
\end{pozn}

\section{A general scheme for path-based models}\label{S3}
On comparing \eqref{DDF-G} and \eqref{HDF-G}, we observe that the DDF-g and HDF-g models are qualitatively similar, differing only in the nature of their output bounds, which are linear for DDF-g but nonlinear for HDF-g. The latter is formalised through the function $\psi^y(\theta)=\frac 1\theta$.

We shall now generalise this approach in the direction of \emph{convex} functions $\psi^y$ on the output side and \emph{concave} functions $\psi^x$ on the input side.  Our new model will encompass all the models mentioned up to this point, allowing for a unified analysis of their properties. At the same time, generalisation will allow the construction of new models whose properties can be specified in advance.

By the general scheme (GS) model applied to $(\x_o,\y_o) \in \T$ with directions $\g_o=(\g_o^x,\g_o^y)\gneqq0$ that may depend on $(\x_o,\y_o)$,
 we understand
\begin{subequations}\label{general}
\begin{align}
(\text{GS})_o\qquad\min\ &{} \theta \label{general1}\\
     &{} \X\la \le \x_o + (\psi^x(\theta)-1) \g_o^x,\label{general2}\\
     &{} \Y\la \ge \y_o + (\psi^y(\theta)-1)\g_o^y,\label{general3}\\
     &{} \ones^T\la =1, \quad      \la \ge \0.\label{general4}
\end{align}
\end{subequations}
Here, the real functions $\psi$, their domains ($\dom$), and their images ($\im$) satisfy the following assumptions.
\begin{enumerate}[label={\rm(A\arabic{*})}, ref={\rm(A\arabic{*})}]
 \item\label{A1}
 $\dom(\psi^x)=(a^x,\infty)$ with $a^x\in  \{-\infty,0\}$ and $\dom(\psi^y)=(a^y,\infty)$ with $a^y\in  \{-\infty,0\}$;
 \item\label{A1bis} $\psi^x$ is smooth, concave, increasing, and $\psi^y$ is smooth, convex,  decreasing;
 \item\label{A2}  $\psi^x(1)= \psi^y(1)=1$;
 \item\label{A3}  $\text{im}(\psi^x)=(b^x,\infty)$ with  $b^x=-\infty$ if $\g^y=\0$ and $b^x\in \{-\infty,0\}$ otherwise; $\im(\psi^y)=(b^y,\infty)$ with $b^y\in \{-\infty,0\}$.
 \end{enumerate}
Due to the monotonicity of $\psi^x$ and $\psi^y$ assumed in \ref{A1bis}, the assumption \ref{A3} means $\lim_{\theta \to a^x_+}\psi^x(\theta)=b^x$ and $\lim_{\theta \to \infty}\psi^x(\theta)=\infty$ and $\lim_{\theta \to a^y_+}\psi^y(\theta)=\infty$ and $\lim_{\theta \to \infty}\psi^y(\theta)=b^y$.

If the vectors $\g_o^x$ and $\g_o^y$ are positive, then the GS model \eqref{general} operates in the full input-output space, and hence it is a
 graph model in the terminology of \cite{fare.al.85}. With $\g_o^x=0$ or $\g_o^y=0$,%
 \footnote{In the case $\g_o^x=0$ or $\g_o^y=0$ the contribution of $\psi^x$ in \eqref{general2} or $\psi^y$ in \eqref{general3} vanishes; in such cases
we take $\dom(\psi^x)=(-\infty,\infty)$ or $\dom(\psi^y)=(-\infty,\infty)$, respectively.}
one obtains output-oriented or input-oriented models, respectively.%
\footnote{In the case of input/output-oriented GS models (where $\g_o^y=\0$ or $\g_o^x=\0$), the GS model can be converted to a DDF-g input/output-oriented model with the same projection as the GS model.  This observation follows from the fact that the functions $\psi^x(\cdot)$ and $\psi^y(\cdot)$ are assumed to be monotonic  
 so that their inverse functions exist and the transformations $\psi^x(\theta)=:a\iff (\psi^x)^{-1}(a)=\theta$ or $\psi^y( \theta)=:b\iff (\psi^y)^{-1}(b)=\theta $ can be applied. 
 In the case of graph models, the inverse transformation allows one to convert the GS model into a form where the input (resp., output) side is linear in $\theta$.}

Denote the right-hand sides of \eqref{general2} and \eqref{general3} by $\p_o^x$ and $\p_o^y$, respectively, that is,
\begin{equation}\label{pathdef}
\p_o^x(\theta):=\x_o + (\psi^x(\theta)-1)\g_o^x; \qquad
\p_o^y(\theta):=\y_o + (\psi^y(\theta)-1)\g_o^y.
\end{equation}
It is easy to see that the models presented in the previous section are special cases of the GS model. These models, taken in conjunction with the directions in Table~\ref{Tab:choices} and the usual assumptions on the positiveness of the data (that is, condition \eqref{xygneqq0} for the BCC models and \eqref{xy>0} for the hyperbolic model), will be called \emph{standard path-based models}. The corresponding parameterisations in Table \ref{Tab:values} indicate that only two choices of the function $\psi^y$ are associated with the standard models: linear, $\psi^y(\theta)=2-\theta$, and hyperbolic, $\psi^y(\theta)=\frac{1}{\theta}$.  Only one form of $\psi^x$ appears: $\psi^x(\theta)=\theta$.
\begin{table}[t!]
\centering
 \begin{tabular}{l c c c c c c}
 \hline\\[-1.8ex]
 Model  & $\g_o^x$ &$\p_o^x(\theta)$ & $\g_o^y$ & $\p_o^y(\theta)$ & $\psi^x$& $\psi^y$ \\[0.8ex]
 \hline\\[-1.8ex]
 BCC-I & $\x_o$ & $\theta \x_o$ & $\0$ & $\y_o$ & $\theta$& $\psi^y$\ \tablefootnote{Since $\g_o^y=0$, the term $\g_o^y \psi^y $ vanishes.}
  \\[0.2ex]
 BCC-O & $\0$ & $\x_o$ & $\y_o$ & $\frac{1}{\theta}\y_o$ & $\psi^x$\ \tablefootnote{Since $\g_o^x=0$, the term $\g_o^x \psi^x $ vanishes.} & $\frac{1}{\theta}$  \\[0.4ex]
 DDF-g & $\g_o^x$ & $\x_o+(\theta-1) \g_o^x$ & $\g_o^y$ & $\y_o+(1-\theta)\g_o^y$ & $\theta$& $2-\theta$  \\[0.4ex]
 DDF& $\x_o$ & $\theta \x_o$ & $\y_o$ & $(2-\theta)\y_o$ & $\theta$ & $2-\theta$\\[0.4ex]
 HDF& $\x_o$ & $\theta \x_o$ & $\y_o$ & $\frac{1}{\theta}\y_o$ & $\theta$& $\frac{1}{\theta}$  \\[0.4ex]
 HDF-g & $\g_o^x$ & $\x_o+(\theta-1) \g_o^x$ & $\g_o^y$ & $\y_o+(\frac{1}{\theta}-1)\g_o^y$ & $\theta$& $\frac{1}{\theta}$  \\[0.8ex]
 \hline
 \end{tabular}
 \caption{Parameterization of the standard path-based models.}
\label{Tab:values}
\end{table}
However, other choices of $\psi$ that satisfy the assumptions \ref{A1}--\ref{A3} are possible, for example $\psi^x(\theta)=\theta^{-p}$, $-1<p<0$, or $\psi^x(\theta)=1+\ln \theta$, and $\psi^y(\theta)=\theta^{-p}$, $p> 0$,   $\psi^y(\theta)=1-\ln \theta$ with domains $(0,\infty)$,
or $\psi^y(\theta)=\e^{1-\theta}$ with domain $(-\infty,\infty)$.  We will be able to fine-tune the model properties by making suitable choices of $\psi^x$ and $\psi^y$ and directions. Note that the choice of functions $\psi^x(\theta)=\theta^{1-p}, \psi^y(\theta)=\theta^{-p}, p\in [0,1] $ and directions $(\g^x,\g^y)=(\x_o,\y_o)>\0$ leads to the so-called generalised distance function introduced by \cite{chavas.cox.99}.

The next subsection shows that the (GS)$_o$ model is well defined as a convex programme, that is, for any $(\x_o,\y_o) \in \T$  the optimal value $\theta^*$ of (GS)$_o$ is finite and the minimum is attained. To this end, we introduce a geometric interpretation of the model.

\subsection{Geometric interpretation}

\begin{figure}
    \centering
     \scalebox{0.07}
       { \includegraphics{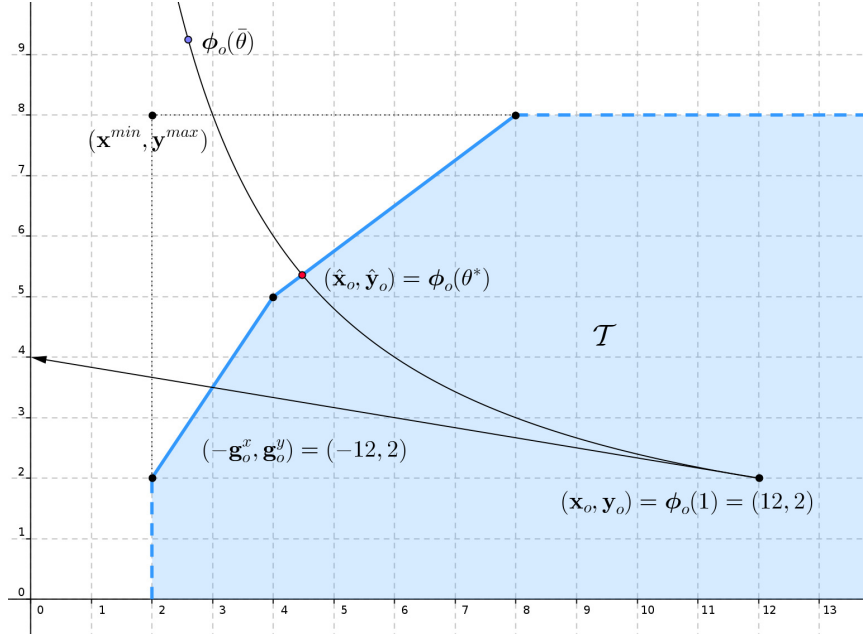}}
    \caption{Hyperbolic path $\p_o(\theta)$ for direction $(12,2)$ as an illustration of  Theorem~\ref{dominance}. The path is a smooth curve that starts for $\theta=1$ at $(\x_o,\y_o)\in\T$ and for decreasing values of $\theta$ it moves in a set of points that gradually dominate each other, until the curve for some $\bar\theta<1$ reaches the value  $\p_o^y(\bar\theta)>y^\M$. As a consequence, the path has to leave $\T$ at some $(\hat\x_o,\hat\y_o)\in\partial \T$  for some $\theta=\theta^*_o$.}
    \label{F:1}
\end{figure}

The definition of the technology set $\T$ in (\ref{T}) and the notation in \eqref{pathdef} allow us to rewrite (GS)$_o$ \eqref{general} in the form
\begin{equation}\label{curve}  \min \{\theta : (\p_o^x(\theta),\p_o^y(\theta))\in {\T}\}.
\end{equation}
In this context, we can interpret the (GS)$_o$ model as a path-based model, where the map
\begin{equation}\label{path}
\p_o: \theta\mapsto (\p_o^x(\theta),\p_o^y(\theta))
\end{equation}
defines a continuous path in the input-output space $\R^m\times \R^s$  parameterised by $\theta\in \cal{D}$, where
\begin{equation}\label{Domain}
    \D=\dom(\psi^x)\cap\dom(\psi^y).
\end{equation}
It is seen from the definition of the path in \eqref{pathdef} and \eqref{path} that each path is determined by a specific choice of $\x_o, \y_o, \g_o^x, \g_o^y$, and $\psi^x$,$\psi^y$ (the \emph{parameters of the path}). Unless otherwise stated, we take \emph{standard parameter choices}, that is,
\begin{equation}\label{eq:spc}
(\x_o, \y_o)\in\T,\quad (g_o^x,g_o^y)\gneqq 0,\quad \text{and}\quad \psi^x, \psi^y \text{ that meet \ref{A1}--\ref{A3}}.
\end{equation}
The next theorem summarises the properties of the path $\phi_o$ that are common for all standard parameter choices. Its proof is placed in  \ref{A}. The path  properties are illustrated  in Figure~\ref {F:1} with an example of a hyperbolic path.

\begin{theorem}\label{dominance}
The path $\p_o$ defined by \eqref{pathdef} and \eqref{path} with standard parameter choices \eqref{eq:spc} has the following properties.
\begin{enumerate}[(a)]
 \item\label{T.a}   For $\theta_1\leq\theta_2\in\D$, the point $\p_o(\theta_1)$ on the path $\p_o$  dominates the point $\p_o(\theta_2)$, (i.e., $\p_o(\theta_1)\succsim \p_o(\theta_2)$),
      and  $\p_o(\theta_1)\neq  \p_o(\theta_2)$, if and only if $\theta_1<\theta_2$.
 \item\label{T.b}  $\p_o(1)=(\x_o,\y_o)\in \T$.
  \item\label{T.c}There exists $i$ such that $\lim_{\theta\to a^x_-}\phi_{io}^x(\theta)=-\infty$ or $r$ such that  $\lim_{\theta\to a^y_+}\phi_{ro}^y(\theta)= +\infty$ (see Figure \ref {F:1} for illustration).
 \item\label{T.d} There exists $\theta_o^*\le 1$ such that $\p_o(\theta ^*)\in \T$ and $\p_o(\theta)\notin \T$ for all $\theta <\theta^*$.
 As a consequence, the point $\p_o(\theta_o^*)$ belongs to the boundary of  $\T$ and  $\theta^*$ is the finite optimal value of the (GS)$_o$ model \eqref{general}.
\end{enumerate}
\end{theorem}

By Theorem~\ref{dominance}, the GS model applied to $(\x_o,\y_o) \in \T$ with directions $(\g_o^x,\g_o^y)\gneqq0$ is well defined.  The optimal value $\theta^*_o$ is called \emph{efficiency score}, or alternatively \emph{value of the efficiency measure} for $(\x_o,\y_o)$. The point $(\p_o^x(\theta_o^*),\p_o^y(\theta_o^*))$ on the path $\p_o$ is called the \emph{projection} of $(\x_o,\y_o)$ in the (GS)$_o$ model.  Hereafter, we shall denote the projection point more compactly by
\begin{equation}\label{projection}(\hat \x_o,\hat \y_o):=(\p_o^x(\theta_o^*),\p_o^y(\theta_o^*)).\end{equation}

Note that although  the projection point $(\hat \x_o,\hat \y_o)\in\partial\T$ is defined uniquely, the uniqueness may not extend to the $\la$-component of an optimal solution $(\la^*,\theta^*)$.
Consequently, the program (GS)$_o$ may have  multiple  optimal slacks given by
$$ (\bm s^{x*}, \bm s^{y*}):=(\hat \x_o-\X\la^*, \Y\la^*-\hat \y_o)\ge \0.
 $$
\begin{lemma}
Let $(\la^*_o,\theta^*_o)$ be an optimal solution of $(GS)_o$. Then at least one of the inequalities in \eqref{general2} or \eqref{general3}  is tight, i.e, at least one component of the corresponding  slacks is zero.
\end{lemma}
\begin{proof}
If all components of the slacks are positive, then due to the assumed continuity of the path there is $\epsilon>0$ such that $(\la^*_o,\theta_o^*-\epsilon)$ is feasible, which contradicts the optimality of  $\theta_o^*$.
\end{proof}

\section{Properties of the general model}\label{S4}

 We will now analyse the (GS)$_o$ model in the light of the ten desirable properties \ref{P1}--\ref{P10}. The only property that the general scheme (GS)$_o$ satisfies regardless of the specific choice of directions is the unique projection property \ref{P1}. This is a simple consequence of Theorem \ref{dominance}. The remaining properties depend on the choice of the direction vector  $g_o$ and the functions $\psi^x$ and $\psi^y$. We examine them one at a time in the following subsections.

\subsection{Indication}
\begin{itemize}
\myitem{(P2)}\label{P2bis} (a) If the evaluated unit is strongly efficient, then the measure equals 1.

(b) If the measure equals 1, then the evaluated unit is strongly efficient.
\end{itemize}

 The  part (a) of this property holds universally in the (GS) scheme as indicated in the following theorem.

 \begin{theorem}\label{weakindication}
   If $(\x_o,\y_o)$ is strongly efficient in ${\T}$, then $\theta^*_o=1$.
    \end{theorem}
 \begin{proof}
   If $(\x_o,\y_o)$ is strongly efficient, then the set of points from $\T$ dominating $ (\x_o, \y_o) $ contains only $ (\x_o, \y_o) $. In view of
  $\p_o(1)=(\x_o,\y_o)$ and Theorem \ref{dominance}\ref{T.a}, one has  $\theta^*_o = 1$.
  \end{proof}
Note that $\theta_o^*=1$ in conjunction with Theorem~ \ref{dominance}\ref{T.d} only yield  $(\x_o,\y_o)\in\partial\T$, hence one must admit the possibility that $(\x_o,\y_o)$ belongs to $\partial^W\T$ in some cases with $\theta_o^*=1$.

\begin{figure}[t]
    \centering
     \scalebox{0.09}
        {\includegraphics{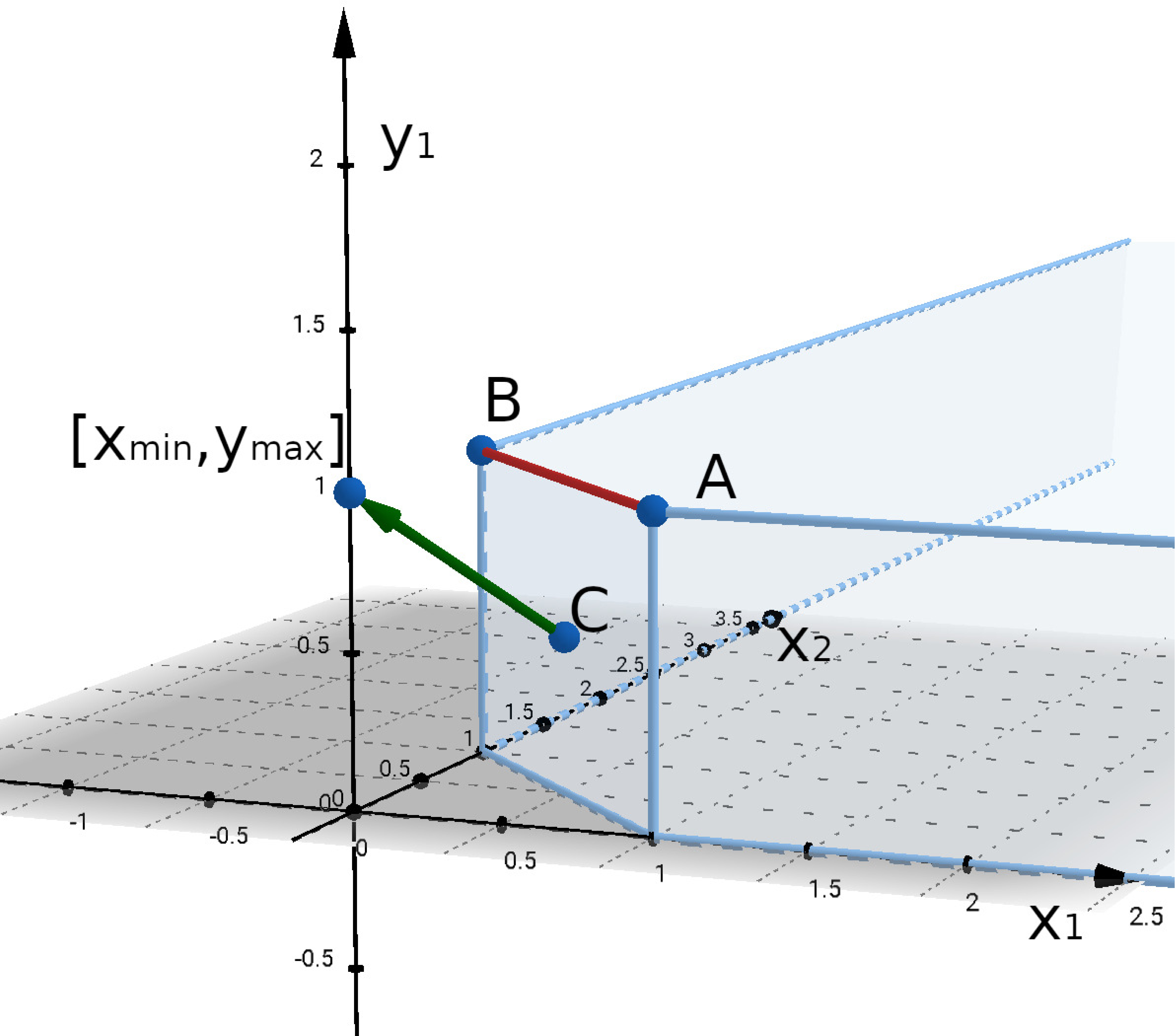}}
    \caption{An example of $\T$, for which the DDF-\ref{G2} model does not project  $C\in \partial^W\T$   onto the strongly efficient frontier $\partial^S\T=\vert AB\vert$. The two-input and single-output technology set is generated by two strongly efficient units, $A=(1,0,1)$ and $B=(0,1,1)$. The unit $C=(0.5, 0.5, 0.5)$ is weakly but not strongly efficient. Here $(\x^\m, \y^\M)=(0, 0, 1)$ and \ref{G2} direction for $C$ is $\g_C=(0.5, 0.5, 0.5)$. Since $\g_C>0$,  Theorem \ref{weakindication1} yields that the part (b) of \ref{P2bis} is not satisfied. The DDF-\ref{G2} path for $C$ with the $\g_C$ direction leaves $\T$ at $C$, which means that $C$  is projected onto itself with $\theta^*_C=1$.}
    \label{F:3}
\end{figure}

\begin{theorem}\label{weakindication1}
The GS model with positive directions violates the part (b) of \ref{P2bis} for any $\T$.
 \end{theorem}
 \begin{proof} Let $(\x_o,\y_o)\in\partial^W\T$ and $\g_o>\0$. We want to show that $\theta_o^*=1$. Assume by contradiction that $\theta_o^*<1$. Then $\hat\x_o=\x_o + (\psi^x(\theta^*_o)-1) \g_o^x<\x_o$ and $\hat\y_o=\y_o + (\psi^y(\theta^*_o)-1) \g_o^y>\y_o$, and $(\hat\x_o,\hat\y_o)\in\T$, which is, by  Lemma \ref{L:weak}\ref{W.c},
 in contradiction with the weak but not strong efficiency of $(\x_o,\y_o)$.
 \end{proof}
\begin{figure}[th]
    \centering
     \scalebox{0.06}
        {\includegraphics{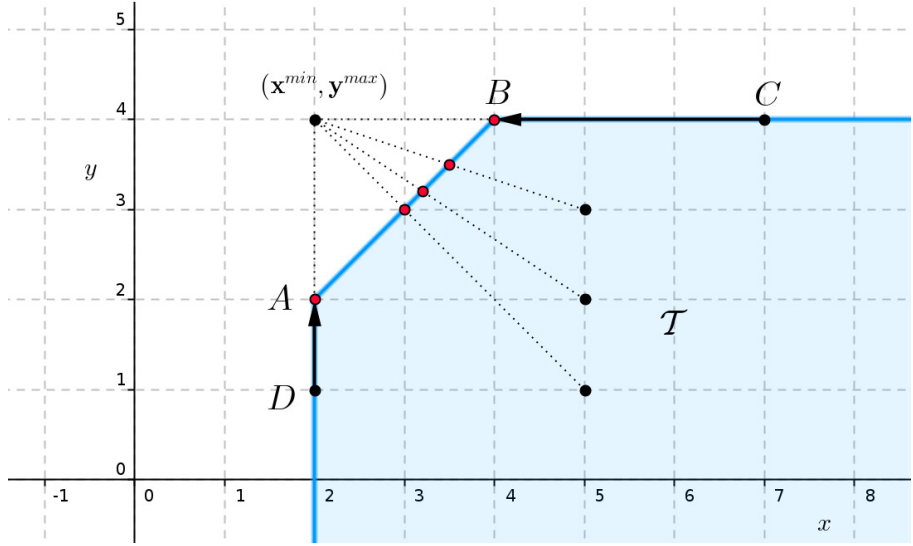}}
    \caption{In  the two dimensional technology set $\T$, the DDF-\ref{G2} model projects units from $\T$ onto $\partial^S\T$. The set  $\T$  is generated by four units. The units $A$ and $B$ are strongly efficient, while $C$ and $D$ belong to $\partial^W\T$. The \ref{G2} direction for $C$ is $(5,0)$ and it points from $C$ to $B$ which implies that $\theta_C^*<1$.  The \ref{G2} direction for $D$ is $(0,1)$ and it points from $D$ to $A$ which implies that $\theta_D^*<1$. Similar observation applies to all units from $\partial^W\T$.}
    \label{F:2}
\end{figure}
 \begin{pozn}\label{poznG}  By Theorem \ref{weakindication1}, all standard path-based models with directions \ref{G1} and \ref{G3}--\ref{G6} violate the  part (b) of \ref{P2bis} for any $\T$. In the case of DDF-\ref{G2}, whether \ref{P2bis} (b) holds depends on the data configuration (the shape of $\T$). It holds in Figure \ref{F:2}, where each unit $(\x_o,\y_o)\in\partial^W\T$ is projected onto $\partial^S\T$ in DDF-\ref{G2} model, and hence $\theta_o^*<1$. A failure of \ref{P2bis} (b) is seen in Figure~\ref{F:3}, where the unit C belongs to $\partial^W\T$ and the corresponding  paths at C point out of $\T$.  This implies that the efficiency scores of C is equal to one. On the other hand,  HDF-\ref{G2} violates  part (b) of \ref{P2bis} even in the case of technology set depicted in Figure \ref{F:2}.  Indeed, since $\frac{1}{\theta}-1>1-\theta$ for $\theta<1$,  the hyperbolic path for $C$ is  above the directional path for $C$  for $\theta<1$, and hence leaves $\T$ at $C$.
 \end{pozn}

The tools to determine whether a unit with an efficiency score equal to one is strongly efficient will be presented in the next subsection, where this problem will be solved in a more general setting.

\subsection{Strong efficiency of projection}
\begin{itemize}
\item[\ref{P3}] The projection point generated by the measure is strongly efficient.
\end{itemize}

The projection in the (GS)$_o$ model is located on the boundary of $\T$ but not necessarily on its strongly efficient part. Note that the lack of strong efficiency of a projection is a known possibility in the BCC and DDF models. In this subsection, we show that this is the case for all standard path-based models.

The absence of strong efficiency of the projection has knock-on effects, such as the lack of strict monotonicity of the efficiency score or the failure of one-to-one identification of efficiency.   Moreover, if a unit is not projected onto the strongly efficient frontier, then its efficiency score does not capture
 all sources of inefficiencies in some input/output components, and hence it overestimates the unit's performance. Therefore, it is important to recognise
whether the specific unit was projected onto a strongly efficient frontier or not. Unfortunately, this recognition cannot be made from the mere value of the optimal score. However, in connection with BCC and DDF,
tools have been developed to detect this feature. Now we are adopting them into our scheme.

   The next theorem states that to identify the strong efficiency of projections, one needs to know more than just the optimal value. The proof is located in \ref{A}.
   \begin{theorem}\label{indication}
   The projection of $(\x_o,\y_o)$ is strongly efficient if and only if for each optimal solution $(\la^*_o,\theta^*_o)$ of \eqref{general}, the inequality constraints \eqref{general2} and \eqref{general3}  are satisfied with equality.
 \end{theorem}

      Theorem \ref{indication} identifies whether the projection is from the strongly efficient frontier theoretically only.
A practical method to identify the strong efficiency of a weakly efficient unit $ (\x_o,\y_o)\in {\T}$ with $\theta^*_o=1$, but also the projection of an inefficient unit $ (\x_o,\y_o)\in {\T}$ with $\theta^*_o< 1$, is called \emph{second phase procedure}. The procedure solves a modified programme that, at a fixed value of $\theta^*_o$, maximises the sum of slacks. In the context of the GS scheme, it reads as follows:
\begin{equation}\label{phase}
\begin{array}{rl}
\max &\ones^T\s^x+\ones^T\s^y  \\
     & \X\la +\s^x =\p_o^x(\theta^*_o),\\
     & \Y\la -\s^y =\p^y_o(\theta^*_o),\\
     & \ones^T\la =1, \
      \la \ge \0, \s^x\ge \0, \s^y\ge\0.
\end{array}
\end{equation}
 \begin{theorem}\label{indication2}
The optimal value in \eqref{phase} vanishes if and only if the projection of $ (\x_o,\y_o)$ in (GS)$_o$ is strongly efficient. Furthermore, for each optimal solution $(\la^*, \s^{x*}, \s^{y*}) $ in \eqref{phase}, the unit $(\X\la^*,\Y\la^*)$ is a strongly efficient benchmark for $ (\x_o,\y_o)$.
\end{theorem}
\begin{proof} The first part follows from Theorem \ref{indication}. The second part claims that $(\X\la^*,\Y\la^*) =(\hat\x_o-\s^{x*},\hat\y_o+\s^{y*})\in\partial^S\T$. Suppose this is not the case. Then by Lemma \ref{L:weak}\ref{W.b}, there exists a unit in $\T$ (corresponding to some $\hat \la\geq \0$ with $\ones^T\hat\la=1$) that dominates but is not equal to $(\X\la^*,\Y\la^*)$. The slacks of $(\hat\x_o,\hat\y_o)$ relative to the new unit have increased, which contradicts the optimality of the original slacks in \eqref{phase}.
\end{proof}

\begin{pozn}\label{benchmark} While in the class of slacks-based models the notions of `projection' and `benchmark' coincide  (see Remark 2 in  \citealp{halicka.trnovska.21}), in the class of path-based models  it is necessary to single out one projection from among multiple benchmarks. The projection is uniquely determined by \eqref{projection} as the point (benchmark) on the path, where the path  leaves  $\partial\T$, which also pins down the optimal value of the objective function. However, the model is capable of generating other benchmarks, i.e., points $(\X\la^*, \Y\la^*)\in\partial T$ that correspond to the alternate optimal solutions $(\la^*,\theta^*_o)$  of (GS)$_o$.  Some of them lie in $\partial^S\T$ and can be generated by optimal solutions of
\eqref{phase}.
\end{pozn}
\begin{pozn}\label{poznG2} Property \ref{P3} implies \ref{P2bis}. Indeed, if $\theta_o^*=1$, then $ (\x_o,\y_o)$ coincides with its projection. Now, if \ref{P3} holds,  this projection is strongly efficient, which yields the (b) part of \ref{P2bis}. The  part (a) of \ref{P2bis} holds universally by Theorem~\ref{weakindication}.

 As a consequence of the implication \ref{P3} $\Rightarrow$ \ref{P2}, the conclusions of Remark \ref{poznG} regarding failure of \ref{P2} for standard path-based models carry over to property \ref{P3}.
\end{pozn}

\subsection{Boundedness}\label{boundedness}
\begin{itemize}
\item[\ref{P4}] The measure takes values between zero and one.
\end{itemize}
If $(\x_o,\y_o)\in {\T}$, then automatically $\theta_o^*\le 1$ (since $\theta=1$ is a feasible solution).  On the other hand, the property $\theta_o^*\ge 0$ is based on the domains of $\psi^x$ and $\psi^y$.  With $\dom(\psi^x)= (0,\infty)$ or  $\dom(\psi^y)= (0,\infty)$,  which only occur for non-linear $\psi^x$ and $\psi^y$, one has $\theta_o^*> 0$ for any choice of $g_o$.
  With linear $\psi$ one has $\dom(\psi)= \R$, and hence \ref{P4} may fail. In the following theorem, we formulate sufficient conditions on $g_o$ for \ref{P4} to hold.
 \begin{theorem}\label{ohranicenost} Consider a GS model with $\psi^x(\theta)=\theta$, $\psi^y(\theta)=2-\theta$, and directions $g_o\gneqq 0$.  If $$\max_i\left\{\frac{x^\m_i+g_{io}^x-x_{io}}{g_{io}^x} : g_{io}^x>0\right \}\ge 0\quad  \text{or}\quad  \max_r\left\{\frac{y_{ro}+g_{ro}^y-y^\M_r}{g_{ro}^y} : g_{ro}^y>0\right\}\ge 0,$$ then $\theta^*_o\ge 0$.
\end{theorem}
\begin{proof} We proceed by constructing lower bounds for the optimal score $\theta^*_o$.  If $g^x_{io}>0$, then the $i$-th component of \eqref{general2}  yields
       $\frac{x^\m_i+g_{io}^x-x_{io}}{g_{io}^x}\le \psi^x(\theta_o^*)=\theta_o^*$.  Similarly,
       if  $g^y_{ro}>0$, then the $r$-th component of \eqref{general3}
      yields $\frac{y^\M_r+g_{ro}^y-y_{ro}}{g_{ro}^y}\ge \psi^y(\theta_o^*)=2-\theta^*_o$
       and hence   $\frac{y_{ro}+g_{ro}^y-y^\M_r}{g_{ro}^y}\le \theta^*_o$. If the maximal lower bound is non-negative, then  $\theta^*_o\geq 0$, which completes the proof.
       \end{proof}
Sufficient conditions for other linear $\psi^x$ and $\psi^y$ follow by an obvious modification of the proof.
\begin{pozn}\label{poznP3} Theorem~\ref{ohranicenost} offers an easy route to verify that DDF-g with directions \ref{G1}--\ref{G3} satisfies boundedness \ref{P4}. A counterexample with negative $\theta^*_o$ for directions \ref{G4}--\ref{G6} is provided by a single-input, single-output model with two units $A=(1,5)$ and $B=(5,1)$. Here, the average direction \ref{G4} is  $(3,3)$ and the standard deviation direction \ref{G5} equals  $(2, 2)$. The scores for the non-efficient unit B are $-1/3$, $-1$, and  $-3$, respectively, for the directions \ref{G4}--\ref{G6}. The HDF-g family (including HDF itself) satisfies \ref{P4} trivially thanks to the non-negative domain of $\psi$.

The results of Theorem \ref{ohranicenost} are in line with the partial results of \cite{sahoo.al.14} related to the input-oriented directional distance model.

 \end{pozn}
\subsection{Unit invariance}\label{UI}
\begin{itemize}
\item[\ref{P5}]
The value of the measure does not depend on the units of measurement in the input and output variables.
\end{itemize}
The next theorem provides a necessary and sufficient condition for unit invariance in terms of invariance of directions with respect to the data transformation. Its proof is given in \ref{A}.
\begin{theorem}\label{invariance}
 Let $\bm C\in \mathbb{D}_{++}^m$ and $\bm B\in \mathbb{D}_{++}^s$ be diagonal matrices with positive elements on the diagonal.  For all $j=1,\dots,n$, consider transformed inputs $\bm C\x_j$ and transformed outputs $\bm B\y_j$, respectively.  Let $\g_o^x$, $\g_o^y$, and $\g^{Cx}_o$, $\g^{By}_o$ denote the directional vectors before and after the data transformation, respectively. The (GS)$_o$ model is unit invariant if and only if 
 one has $\g^{Cx}_o=\bm C\g^{x}_o$ and $g^{By}_o=\bm B\g^{y}_o$.
\end{theorem}

\begin{pozn}\label{poznP5}Theorem~\ref{invariance}  yields easily that GS models with \ref{G1}--\ref{G5} directions satisfy unit invariance \ref{P5} but this property fails with directions \ref{G6}.
\end{pozn}
\subsection{Translation invariance}
\begin{itemize}
\item[\ref{P6}] The value of the measure is not affected by translation in inputs or outputs.
\end{itemize}
{The next theorem provides a necessary and sufficient condition for translation invariance in terms of translation invariance of directions. Its proof is given in \ref{A}.
\begin{theorem}\label{translation}
The (GS)$_o$ model is translation invariant if and only if the directional vectors $\g_o^x$ and $\g_o^y$ are invariant to the translation of inputs and outputs, respectively.
\end{theorem}

The translation invariance of the DDF-g models is analysed in detail by \cite{aparicio.al.16}. This subsection is consistent with their results and extends them to the whole class of GS models.
\begin{pozn}\label{poznP6}Using Theorem~\ref{translation}, one easily verifies that GS models with directions \ref{G2}, \ref{G3}, \ref{G5}, and \ref{G6} satisfy translation invariance \ref{P6} but those with directions \ref{G1} and \ref{G4} do not.
\end{pozn}

\subsection{Monotonicity}\label{ssmon}
\begin{itemize}
\item[\ref{P7}] An increase in any input or a decrease in any output relative to the evaluated unit, keeping other inputs and outputs constant, reduces or maintains the value of the measure.
\end{itemize}
For any fixed $\bar{\theta}\in\cal D$, we consider a function $\p(\bar{\theta})$: $(\x_o,\y_o)\mapsto \p_o(\bar{\theta})$ mapping each point $(\x_o,\y_o)\in \T$ to the point
 $(\p_o^x(\bar{\theta}), (\p_o^y(\bar{\theta}))$  on the path $\p_o(\theta)$ at $\theta=\bar{\theta}$.
The monotonicity property of the models states that for any $(\x_o,\y_o)$ and $(\x_q,\y_q)$ in $\T$ one has  
$$
(\x_o,\y_o)\succsim (\x_q,\y_q)\ \Rightarrow\ \theta_o^*\ge \theta_q^*.
$$
The monotonicity of a model is linked with the monotonicity of $\p(\theta)$, which we define as follows.
\begin{definition}
We say that $\p(\bar{\theta})$ is monotone on $\T$ in $\bar{\theta}\in \mathcal{D}$ if for any two units $(\x_o,\y_o)$ and $(\x_q,\y_q)$ in $\T$, one has
\begin{equation}\label{monoton}
    (\x_o,\y_o)\succsim (\x_q,\y_q) \ \Rightarrow
    \  \p_o(\bar{\theta})\succsim \p_q(\bar{\theta}).
\end{equation}
\end{definition}
For a function $\p(\bar{\theta})$ such that $[\p^x(\bar{\theta})]_i$
depends only on $\x_{io}$ and $[\p^y(\bar{\theta})]_r$ depends only on 
$\y_{ro}$, the 
monotonicity property in Definition~\ref{monoton}
simply means that $[\p^x(\bar{\theta})]_i$ and $[\p^y(\bar{\theta})]_r$ are nondecreasing
in $\x_{io}$ and $\y_{ro}$, respectively. 

\begin{lemma}\label{monotonlema}
Let $(\x_o,\y_o)$ and $(\x_q,\y_q)$
be two units in $\T$ with the corresponding optimal values $\theta^*_o, \theta^*_q$. If $\p_o(\theta^*_o)\succsim \p_q(\theta^*_o)$, then each optimal solution
$(\theta_o^*,\la_o^*)$  of (GS)$_o$ is a feasible solution of (GS)$_q$. Furthermore, $\theta^*_q\le \theta^*_o$.
\end{lemma}
\begin{proof}
The assumptions of the lemma and feasibility of $(\theta_o^*,\la_o^*)$  for (GS)$_o$ yield
$$
\X\la_o^*\le\p_o^x(\theta^*_o)\le \p_q^x(\theta^*_o), \quad \Y\la_o^*\ge \p_o^y(\theta^*_o)\ge \p_q^y(\theta^*_o).
$$
Therefore $(\theta_o^*,\la_o^*)$ is also feasible for (GS)$_q$. 
Since $\theta^*_q$ is the optimal value for the minimisation problem (GS)$_q$, we get $\theta^*_q\le \theta^*_o$.
\end{proof}
The next theorem is a consequence of Lemma \ref{monotonlema}. Here, we assume that the assumption of monotonicity is satisfied only at $\bar\theta\in\tilde{\cal D}$, where $\tilde{\cal D}$ is a set of those $\bar\theta$ for which there exists
$(\x_o,\y_o)\in\T$ such that $\theta_o^*=\bar\theta$.
\begin{theorem}\label{monotonicity}
Suppose that $\p({\bar\theta})$
is monotone at any ${\bar\theta}\in \tilde{\cal D}$.  Then the GS model satisfies the property of monotonicity \ref{P7}.
\end{theorem}
\begin{pozn}\label{poznP7} Note that if  $\g_o$
does not depend on $(\x_o,\y_o)$, then $\p(\theta)$ is monotone for any choice of $\psi^x$ and $\psi^y$ in any $\theta\in \mathcal{D}$.
Therefore, by Theorem \ref{monotonicity}, (GS) models with directions \ref{G3}--\ref{G6} satisfy the property \ref{P7}. For the \ref{G1} direction with arbitrary (standard) $\psi^x$ and $\psi^y$, the monotonicity assumption of Theorem \ref{monotonicity} is satisfied in $\bar\theta\in \mathcal{D}$, provided that $\psi^x(\bar\theta)>0$ and $\psi^y(\bar\theta)>0$ hold. These conditions are satisfied by the DDF-\ref{G1} and HDF-\ref{G1} models, since, by Remark \ref{poznP3}, in these models $\theta^*$ reaches only positive values on $\T$.  Therefore, both models meet the property \ref{P7}.
 For the \ref{G2} direction with arbitrary (standard) $\psi^x$ and $\psi^y$, the monotonicity holds at $\bar\theta\in\cal D$ provided that $\psi^x(\bar\theta)>0$ and $\psi^y(\bar\theta)<2$ hold.  These conditions are again satisfied by both DDF-\ref{G2} and HDF-\ref{G2}.  In fact, according to Remark \ref{poznP3}, DDF-\ref{G2} meets \ref{P4}, and thus $\theta_o^*>0$. Furthermore, $\theta^*_o\ge \frac{1}{2}$ on $\T$ in HDF-\ref{G2} with easy computation.

In conclusion, all standard path-based models have the property of monotonicity \ref{P7}.
\end{pozn}

\subsection{Strict monotonicity}

\begin{itemize}
\item[\ref{P8}] An increase in any input or a decrease in any output relative to the evaluated unit, keeping other inputs as well as outputs constant, reduces the value of the measure.
\end{itemize}
We begin with a necessary condition for strict monotonicity \ref{P8}.

\begin{theorem}\label{strongmonotonicity}
If there exists $(\x_o,\y_o)\in \partial^W\T$ such that $\theta^*_o=1$, then the GS model does not meet the property of strict monotonicity.
\end{theorem}
\begin{proof} Since $(\x_o, \y_o)$ is not a strongly efficient unit, by Theorem \ref{indication2} there exists an optimal solution $(\la^*, \bm s^{x*},\bm s^{y*})$, of \eqref{phase}  such that $(\bm s^{x*},\bm s^{y*})\neq \0$ and $(\X\la^*,\Y\la^*)\in \partial^S\T$. Since $(\x_p, \y_p)=(\X\la^*,\Y\la^*)$ is strongly efficient,  Theorem \ref{weakindication} yields $\theta^*_p=1$. Strict monotonicity is violated because $(\x_p,\y_p)$ dominates and is different from $(\x_o, \y_o)$ but $\theta^*_p=\theta_o^*=1$.
\end{proof}
We have shown in Theorem~\ref{weakindication1} that positive directions rule out the validity of \ref{P2}. As a consequence of Theorem \ref{strongmonotonicity}, an analogous statement applies to \ref{P8}.
\begin{theorem}\label{T:positiveP8}
The GS model with positive directions fails in \ref{P8} for any $\T$ defined by \eqref{T}.
\end{theorem}
\begin{pozn}\label{poznP8}
According to Remark \ref{poznG},  
all  standard path-based models admit $(\x_o,\y_o)\in \partial^W\T$ such that $\theta^*_o=1$ for some data configurations and hence by Theorem \ref{strongmonotonicity} do not satisfy \ref{P8}.
\end{pozn}

 \subsection{Super-efficiency}

\begin{itemize}
\item[\ref{P9}] The value of the measure of a unit outside the technology set is well defined and finite.
\end{itemize}
It is well understood that for oriented models super-efficiency may fail for \emph{some} $(\x_o,\y_o)$. Specifically, if $(\x_o,\y_o)\notin \T$, some of the inequalities in \eqref{general2}  or \eqref{general3} are not satisfied for $\theta=1$. To remedy the situation, one must increase the value of $\theta=1$ and thus increase $\phi_o^x(\theta)$ on the input side and decrease $\phi_o^y(\theta)$ on the output side. This increase/decrease is impossible if any of the components of $\phi_o^x(\theta)$ or $\phi_o^y(\theta)$ are not dependent on $\theta$. It follows that if one wishes to ensure \ref{P9} for arbitrary $(\x_o,\y_o)\notin \T$, one must consider only the graph models where
$\g_o^x>0$ and $\g_o^y>0$. The following theorem establishes  some necessary and some sufficient conditions  for \ref{P9}. Its proof is located in \ref{A}.

\begin{theorem}\label{superefficiency}
For $(\x_o,\y_o)\notin \T$ and $(\g_o^x,\g_o^y)\gneqq0$, the following statements hold.
\begin{enumerate}[(a)]
\item\label{superefficiency.a} If (GS)$_o$ is feasible, that is, if there exists $\bar\theta\in\dom(\psi)$  such that $(\p_o^x(\bar\theta),\p_o^y(\bar\theta))\in\T $,
then (GS)$_o$  admits an optimal
solution and the optimal value satisfies $ \theta^*_{o}> 1$.
\item\label{superefficiency.b}
 Let $\g_{o}>0$. If
 $\y^{\min}-\y_o+\g_o^y>0$ or
 $\im(\psi^y)=(-\infty, \infty)$,
 then (GS)$_o$  admits an optimal
solution and the optimal value satisfies $ \theta^*_{o}> 1$.
 \item\label{superefficiency.c}
 Assume  $\im(\psi^y)=(0,\infty)$ and let $\g_o> 0$.   If (i) there is no $\la\ge 0$ 
 with $\ones^T\la =1$ such that inequality $\Y\la\ge \y_o-\g_o^y$ is satisfied or (ii) ${y^\M_r-y_{ro}+g_{ro}^y}<0$
for some $r$,
 then (GS)$_o$  is infeasible. 
 \end{enumerate}
\end{theorem}

\begin{pozn}
With $\g_o>0$, DDF-g models are feasible for all units outside of the technology set by Theorem \ref{superefficiency}\ref{superefficiency.b} ($\dom(\psi)=(-\infty, \infty)$). As a result, the entire family of DDF-g models satisfies the property \ref{P9} provided that the directional vectors are positive.  Note that the positivity of the directions \ref{G2} for $(\x_o,\y_o)\notin \T$ is satisfied only if $\x_o> \x^{\m}$ and $\y_o< \y^{\M}$.

In the class of HDF-g models over positive data, only the HDF model with positive \ref{G1} directions is super-efficient since it satisfies the assumptions of Theorem \ref{superefficiency}\ref{superefficiency.b}.

The HDF-g models over positive data with \ref{G3}--\ref{G6} directions do not guarantee \ref{P9}. In fact, one can always choose $(\x_o,\y_o)\notin\T$ with $y_{ro}>y^{\M}_r+g_{ro}$ for some $r$.  The infeasibility of the corresponding $(GS)_o$ then follows from (ii) of Theorem \ref{superefficiency}\ref{superefficiency.c}.

Condition (ii) is too coarse for
HDF-\ref{G2} because inequality $y_{ro}>y^{\M}_r+g_{ro}$ is incompatible with positive \ref{G2} directions. Instead, we shall provide an example that meets condition (i) of Theorem \ref{superefficiency}\ref{superefficiency.c}. Consider a one-input, two-output production technology with the technology set generated by two units $A=(1,1,10)$ and $B=(1,10,1)$. It can be easily shown that the unit $(\x_o,\y_o)=(2,8,8)$ is outside $\T$. The corresponding \ref{G2} direction is equal to $\g_o=(1, 2, 2)$. An easy calculation now shows that condition (i) of Theorem \ref{superefficiency}\ref{superefficiency.c} is satisfied, and hence HDF-g with \ref{G2} directions does not guarantee super-efficiency \ref{P9}.
\end{pozn}
\begin{pozn}
The definition of the technology set $\T$ in \eqref{T}, even when applied to non-negative data, does not presuppose non-negative values of outputs. Consequently, for the purposes of super-efficiency measurement, we allow projections of units $0\leq(\x_o,\y_o)\notin \T$ onto elements of the frontier with negative outputs.
Therefore, our results differ from those of studies by other authors, where projections onto the negative part of the frontier are treated by definition as \emph{infeasible}.

For example, in \cite{briec.kerstens.09}  (in the context of general economic productivity theories), DDF-g is classified as not having property \ref{P9} only because some units project onto the negative part of the frontier. To rule out such a source of infeasibility, the authors separately consider models, where the technological set is extended by the so-called \emph{free disposal cone}, which then matches our approach.   \cite{johnson.mcginnis.09} praise the HDF model, which in contrast to the DDF-g models is always feasible for positive data.
This can also be seen in our GS scheme, since the HDF projection with \ref{G1} directions reads $(\theta \x_o, \frac{1}{\theta}\y_o)$ and thus is positive for positive data. \cite{mehdiloozad.roshdi.14.arxiv} provide a comprehensive review of super-efficient DEA models and analyse the pivotal role of the directional vector for super-efficiency.
\end{pozn}
\subsection{Homogeneity}\label{s:homogeneity}
\begin{itemize}
\item[\ref{P10}]  {\bf Homogeneity.} Feasible scaling of the input vector and the output vector of the evaluated unit by a power of $\mu>0$ (the power may be different for inputs and outputs) results in the efficiency score being scaled by a power of $\mu$.
\end{itemize}

In this section, the efficiency score  
provided by the GS model for a point $(x,y)\in \T$ 
will be denoted by $\theta^*(\x,\y)$. Here, the argument in $\theta^*$ indicates that the efficiency score depends on $(\x,\y)$  not only directly, but also indirectly through unit-dependent directional vectors $\g$.

The next definition is consistent with the definition of the so-called almost homogeneity used in connection with the measure of hyperbolic efficiency in \cite{cuesta.zofio.05}. 
\begin{definition}\label{defH}
 We say that the measure $\theta^*$ generated by a GS model over a technology set $\T$ is homogeneous of degree $(\alpha,\beta,\gamma)$ if the relation
\begin{equation}\label{homog}
 \theta^*(\mu^{\alpha}\x, \mu^{\beta}\y)=\mu^{\gamma}\theta^*(\x,\y)
\end{equation}
holds for each $(\x,\y)\in \T$ and for each $\mu>0$ such that the value $\theta^*(\mu^{\alpha}\x, \mu^{\beta}\y)$ is well defined.
\end{definition}
 Note that a measure is homogeneous of degree $(\alpha,\beta,\gamma)$, if and only if it is also homogeneous of degree $(k\alpha,k\beta,k\gamma)$ for any $k\ne 0$. To see this, 
 it suffices to set $\bar\mu=\mu^k$ in \eqref{homog} (the map $\mu\to \mu^k$ is one-to-one for $\mu>0$).  Therefore, without loss of generality, we consider only the homogeneity of degree $(\alpha,\beta, 1)$.
\begin{pozn}
By Definition \ref{defH}, $(\x,\y)\in \T$, and therefore $\theta^*(\x,\y)$ is well defined (Theorem \ref{dominance}\ref{T.d}), that is, the minimum in \eqref{general} is attained as a finite value.
However, the scaled point $(\mu^{\alpha}\x, \mu^{\beta}\y)$ is not necessarily in $\T$. If the programme \eqref{general} for $(\mu^{\alpha}\x, \mu^{\beta}\y)$  is not feasible,
then $\theta^*(\mu^{\alpha}\x, \mu^{\beta}\y)=+\infty$, and \eqref{homog} is not satisfied. On the other hand, if the programme \eqref{general} for $(\mu^{\alpha}\x, \mu^{\beta}\y)$ is feasible, then $\theta^*(\mu^{\alpha}\x, \mu^{\beta}\y)$ is well defined (Theorem \ref{superefficiency}\ref{superefficiency.a}).
\end{pozn}

For the purposes of analysing the homogeneity property, we will restrict ourselves to the technology set $\T_+=\T\cap (\mathbb{R}^m_+\times \mathbb{R}^s_+)$,  generated with positive data. We will say that the GS model is homogeneous of degree $(\alpha, \beta,1)$ if the corresponding measures $\theta^*$ are homogeneous of degree $(\alpha, \beta,1)$ for any $\T_+$.

Next, we show that in the class of the GS models,  the homogeneity property is satisfied exclusively for models, where the functions $\psi^x$ and $\psi^y$ are of specific forms, and the directions are \ref{G1} directions, in the graph models, or $(\x_o,\0)$ in the input-oriented and $(\0,\y_o)$ in the output-oriented models. 

We start with the following theorem that provides a necessary and sufficient condition for homogeneity of the GS model with \ref{G1} directions. The proof can be found in \ref{A}. 

\begin{theorem}\label{homogT} 
The GS model with \ref{G1} directions is homogeneous of degree $(\alpha,\beta,1)$ if and only if 
\begin{equation}\label{homogpath}
    \psi^x(\theta)=\theta^{-\alpha}, \ \alpha\in [-1,0), \qquad
     \psi^y(\theta)=\theta^{-\beta}, \ \beta>0, \quad \text{with} \quad \dom(\psi^x)=\dom(\psi^y)=(0,\infty).
\end{equation}
\end{theorem}

The next theorem provides a necessary and sufficient condition for the homogeneity of input- and output-oriented models. The proof is analogous to the proof of Theorem~\ref{homogT} and hence is omitted. 

\begin{theorem}\label{homogoriented}
\begin{enumerate}[(a)]
\item The input-oriented GS model with $(\g^x,\g^y)=(\x,\0)$ is homogeneous of degree $(\alpha,0,1)$ if and only if 
$   \psi^x(\theta)=\theta^{-\alpha},  \theta>0, \ \alpha\in [-1,0)$.
\item The output-oriented GS model with $(\g^x,\g^y)=(0,\y)$ is homogeneous of degree $(0,\beta,1)$ if and only if $\psi^y(\theta)=\theta^{-\beta}, \ \theta>0,\ \beta>0$.
\end{enumerate}
\end{theorem}

Finally, the following theorem states that models with the \ref{G2}-\ref{G6} directions discussed in this article are not homogeneous. Its proof is in \ref{A}.
\begin{theorem}\label{nonhomog}
The GS model with directions $(\g_o^x,\g_o^y)\gneqq0$ of the form $(\g^x,\g^y)=(\mathbf{h}^x,\mathbf{h}^y)$ or $(\g^x,\g^y)=(\gamma \x_o-\mathbf{h}^x,\mathbf{h}^y-\delta \y_o)$, where 
$(\mathbf{h}^x,\mathbf{h}^y)$ is a constant vector,  
is not homogeneous.
\end{theorem}

\begin{pozn}\label{homhyp}
Theorems \ref{homogT} and \ref{homogoriented} cover the known results for the HDF, BCC-I, and BCC-O models. Necessarily, in these cases one has $\psi^x(\theta)=\theta$ and $\psi^y(\theta)=\theta^{-1}$. Thus,  $\alpha =-1$ and $\beta=1$.  This implies that HDF, BCC-I and BCC-O are homogeneous of degrees $(-1,1,1)$, $(-1,0,1)$, and $(0,1,1)$, respectively.  
These results were originally shown in \cite{fare.al.85} and later in  \cite{cuesta.zofio.05}.
\\
The measure defined by $\psi^x(\theta)=\theta, \psi^y(\theta)=\theta^{-p}$, where $p>0$, is homogeneous of degree $(-1,p,1)$. Furthermore, the measure of the generalised distance function, defined by $\psi^x(\theta)=\theta^{1-p}, \psi^y(\theta)=\theta^{-p}, p\in [0,1]$ and
introduced in \cite{chavas.cox.99}
is homogeneous of degree $(p-1,p,1)$, see also Proposition 2 in \cite{chavas.cox.99}. 
\end{pozn}

\begin{pozn}[g-homogeneity]
 The standard formulation of DDF-g models is in terms of an inefficiency measure $\delta$. This measure can be obtained from our GS models by the transformation $\delta=1-\theta $.  It is well known that  $\delta$ satisfies the so-called \emph{g-homogeneity} (\cite{hudgins.primont.07}) of degree $(1,-1)$ with the following meaning: 
 scaling each direction $\g$ by $\mu>0$ results in scaling the measure $\delta$ by $\mu^{-1}$. 
 
 Let us note that the concept of g-homogeneity of degree $(1,-1)$ could be generalised to more general powers of $\mu$, i.e. to g-homogeneity of degree $(p,q)$. However, a similar procedure as used in the proof of Theorem \ref{homogT} shows that the GS model is g-homogeneous of degree $(p,q)$ if and only if $p=1, q=-1$ and $\psi^x=\theta$, $\psi^y=2-\theta$, which corresponds to linear DDF-g models.
 \end{pozn}

\section{Extensions}\label{S5}

Section~\ref{S4} has introduced some theoretical tools useful in analysing the ten desirable properties \ref{P1}--\ref{P10}. With the exception of \ref{P10}, the properties were analysed over arbitrary data. In the remarks of Section~\ref{S4}, these tools were applied to standard path-based models over non-negative data. Next, we expand the analysis to examine the first nine  properties of standard models on arbitrary data (Section~\ref{SS5.1}).
Furthermore, we will demonstrate how the theoretical tools of Section~\ref{S4} 
facilitate the design of new models, where all or almost all desirable properties are satisfied (Subsection~\ref{SS5.2}).

\subsection{Properties of standard path-based models over arbitrary data}\label{SS5.1}

The particular choices of the vectors $\g_o$ in Table~\ref{Tab:properties} are based on the list found in Table~\ref{Tab:choices}. Since we now allow for negative data, the non-negativity  of directions is ensured by taking absolute values where necessary. This is the case of expressions for \ref{G1} a \ref{G4}, but also for \ref{G2}, where the values of units outside $\T$ could harm the assumption $\g_o\gneqq \0$. Modified directions with absolute values are indicated by \ref{|G1|}, \ref{|G2|}, and \ref{|G4|}.

\begin{table}[t!]
\centering
 \begin{tabular}{l c c c c c c}
  \hline\\[-1.8ex]
  Direction notation & \newtag{\textbar G1\textbar}{|G1|}&\newtag{\textbar G2\textbar}{|G2|}&(G3)&\newtag{\textbar G4\textbar}{|G4|}&(G5)&(G6)\\
 \hspace{4.8em} $\g_o^x$ & $\vert \x_o\vert$  &  $\vert \x_o- \x^{\m}\vert$ &  $  \x^{\M}- \x^{\m}$  &  $ \vert \x^{ev}\vert$ & $ \x^{sd}$ &$\ones$\\[0.4ex]
\hspace{53pt}   $\g_o^y$ &   $\vert \y_o\vert$ &   $ \vert \y^{\M}-\y_o\vert$ &   $ \y^{\M}- \y^{\m}$ &  $\vert \y^{ev}\vert$ & $ \y^{sd}$ &$\ones$\\[0.8ex]
 \hline\\[-1.8ex]
 \ref{P2} {indentification}& \xmark(\xmark)& \xmark$^*$(\xmark)& \xmark(\xmark)& \xmark(\xmark)& \xmark(\xmark)& \xmark(\xmark)\\
 \ref{P3}  strong efficiency  & \xmark(\xmark) & \xmark$^*$(\xmark) & \xmark(\xmark) &  \xmark(\xmark)  &  \xmark(\xmark)&\xmark(\xmark)\\
\ref{P4} $\theta^*_o\in [0,1]$ & \checkmark$^{**}$(\checkmark) & \checkmark(\checkmark) & \checkmark(\checkmark)& \xmark(\checkmark)  & \xmark(\checkmark)&\xmark(\checkmark)\\
\ref{P5} unit invariance
  & \checkmark(\checkmark) & \checkmark(\checkmark) & \checkmark(\checkmark) & \checkmark(\checkmark) & \checkmark(\checkmark)&\xmark(\xmark) \\
\ref{P6} translation invariance& \xmark(\xmark) & \checkmark(\checkmark) & \checkmark(\checkmark) & \xmark(\xmark) & \checkmark(\checkmark)& \checkmark(\checkmark)\\
\ref{P7} {monotonicity} & \checkmark$^{**}$(\checkmark)&\checkmark(\checkmark)&\checkmark(\checkmark)&\checkmark(\checkmark)&\checkmark(\checkmark)&\checkmark(\checkmark)\\
\ref{P8} {strict monotonicity}& \xmark(\xmark)& \xmark$^*$(\xmark)& \xmark(\xmark)& \xmark(\xmark)& \xmark(\xmark)& \xmark(\xmark)\\
\ref{P9} {super-efficiency}&\checkmark$^{*}$(\checkmark$^{**}$)&\checkmark$^{*}$(\xmark)&\checkmark$^{*}$(\xmark)&\checkmark$^{*}$(\xmark)&\checkmark$^{*}$(\xmark)&\checkmark(\xmark)\\
\ref{P10} homogeneity&\xmark(\checkmark$^{**}$)&\xmark(\xmark)&\xmark(\xmark)&\xmark(\xmark)&\xmark(\xmark)&\xmark(\xmark)\\[0.8ex]
 \hline
 \end{tabular}
 \caption{Properties of the DDF-g model for $\psi^x(\theta)=\theta$, $\psi^y(\theta)=2-\theta$  and  the HDF-g model for $\psi^x(\theta)=\theta$, $\psi^y(\theta)=\frac{1}{\theta}$ (in brackets) with respect to different choices of vector pairs $\g_o^x,\g_o^y$.
\checkmark$^{*}$ -- the property is satisfied for positive directions; \checkmark$^{**}$ -- the property is satisfied for positive data but not for general data; \xmark$^*$ -- the property holds only for some very specific $\T$.}
\label{Tab:properties}
\end{table}

The list of desirable properties in Table~\ref{Tab:properties} omits the property \ref{P1}, which in the (GS) scheme is satisfied regardless of the choice of $\psi^x$, $\psi^y$ or $\g_o$. Observe that all the models in Table~\ref{Tab:properties} together fail on the same properties \ref{P2}, \ref{P3}, and \ref{P8}.%
\footnote{It may happen that  DDF-\ref{G2} meets \ref{P2}, \ref{P3}, and \ref{P8} over some very specific $\T$ --- see Figure~ \ref{F:2}.}
Since \ref{P2} is not satisfied in any of the models considered, failure of \ref{P3} is inevitable in view of the implication \ref{P3} $\Rightarrow$ \ref{P2} shown in Remark~\ref{poznG2}.
\begin{itemize}[leftmargin=*]
 \item  {\bf Strong efficiency of projections \ref{P3}.} According to Remark~\ref{poznG2}, all models fail \ref{P3} with non-negative data and therefore also for arbitrary data.

 \item \textbf{Boundedness}  \ref{P4}:  For $\psi^x(\theta)=\theta$, $\psi^y(\theta)=\frac{1}{\theta}$, allowing arbitrary data does not produce a change; for all directions $\g_o\gneqq 0$ one again obtains $\theta^*_o>0$. Similarly, for $\psi^x(\theta)=\theta$, $\psi^y(\theta)=2-{\theta}$, and directions \ref{|G2|}=\ref{G2} and \ref{G3} we have $\theta^*_o\geq0$ by  
 Theorem~\ref{ohranicenost}. However, for directions \ref{|G1|}, the status changes; negative data may lead to a negative score, as shown in the example of one input, one output with two units $A=(-1,3), B=(3,1)$, where the directional vector for $B$ is $(3,1)$ and its score equals $-1/3$. However, if at least one of the inputs takes only positive values, then Theorem~\ref{ohranicenost} implies that $\theta^*_o>0$ for all $(\x_o,\y_o)\in\T$. This is in agreement with the results of \cite{kerstens.vandewoestyne.11}.

 By Remark~ \ref{poznP3}, DDF-g models with directions \ref{P4}--\ref{P6} do not guarantee non-negative efficiency scores for positive data and, therefore, the same is true for general data.

 \item {\bf Unit and translation invariance \ref{P5} and \ref{P6}.} There is no change compared to the positive data case and the conclusions of Remarks \ref{poznP5},  \ref{poznP6} remain valid.

\item {\bf Monotonicity \ref{P7}.}  Allowing negative data does not affect the status of \ref{P7} for directions that do not explicitly depend on $(\x_o,\y_o)$. Therefore, the models with directions \ref{G3}, \ref{|G4|}, \ref{G5}, and \ref{G6} satisfy \ref{P7}. Likewise, the directions
 \ref{|G2|}=\ref{G2} meet the assumptions of Theorem~\ref{monotonicity}, and therefore here there is also no change.

The situation is different for DDF with the directions \ref{|G1|}, which do not satisfy the assumptions of Theorem~\ref{monotonicity}. To see that \ref{P7} may fail, consider a one-input, one-output  counterexample with three units: efficient unit $A=(-1,\ 4)$ and two inefficient units $B=(0.1, \ 0)$, $C=(1,\ 0)$. Although $B$ dominates $C$, one has $\theta^*_B=-10<\theta^*_C=-1$.
We remark that directions $\vert \x_o\vert$, resp. $\vert\y_o\vert$ appear in \cite{cheng.al.13} in connection with radial oriented models. There, it is observed that monotonicity may fail if zero is an internal point of the data.

 \item {\bf Strict monotonicity \ref{P8}.}
 According to Remark~\ref{poznP8}, all models fail \ref{P8} with non-negative data and therefore also for arbitrary data.

  \item {\bf Super-efficiency \ref{P9}.} DDF-g models with directions \ref{|G1|}, \ref{|G2|}, \ref{G3}, \ref{|G4|}, \ref{G5}, and \ref{G6} retain the property \ref{P9} for arbitrary data (as long as the directions are positive). Similarly, the negative results for HDF-g models with the directions \ref{|G2|}, \ref{G3}, \ref{|G4|}, \ref{G5}, and \ref{G6} carry over to arbitrary data. Theorem \ref{superefficiency} does not resolve the status of HDF-g models with \ref{|G1|} directions. However, one can find an example where the model is infeasible. Consider a single-input, single-output case, where technology $\T$ is generated by a single unit $A=(2,0)$. For  $B=(1,10)\notin\T$, whose  \ref{|G1|}  direction reads  $\g_B=(1,10)$, the output inequality $0\ge 10\,\theta^{-1} $ is not satisfied for any $\theta>0$.
  \item {\bf Homogeneity \ref{P10}.} Theoretical tools are derived only for positive data in Section \ref{s:homogeneity}. Nevertheless, it can be shown that the statement of Theorem \ref{nonhomog} can be extended to all GS models over negative data.
 \end{itemize}

\subsection {Non-standard directions with good properties}\label{SS5.2}

 We will focus our attention on the \ref{G2} directions, which we have already shown produce the invariance properties \ref{P5} and \ref{P6} for arbitrary $\psi^x$ and $\psi^y$ satisfying \ref{A1}--\ref{A3} (see Remarks \ref{poznP5} and \ref{poznP6}). These directions allow for a small modification that will also ensure
   the validity of \ref{P5} and \ref{P7}. The idea is to choose $\g_o$ so that the corresponding path $\p_o$ crossing $(\x_o,\y_o)$ at $\theta=1$ passes through the point $(\x^\m,\y^\M)$ at $\theta=  \theta_{\m}$ (see Figure~ \ref{F:3}). Then $\theta_{\m}$ serves as the lower bound for the efficiency score, as indicated in the next lemma. Its proof follows from a simple calculation and is therefore omitted.
\begin{lemma}\label{dir} Let $\theta_{\m}\in [0,1)\cap\mathcal{D} $ and $(\x_o, \y_o)\in\T\setminus\{ (\x^\m, \y^\M)\}$. The path $\p_o(\theta)$ runs through $(\x^\m, \y^\M)$ at $\theta=\theta_{\m}$, i.e., $\p_o({\theta}_{\m})=(\x^\m,\y^\M)$
if and only if
the directions $\g_o^x, \g_o^y$ satisfy
 \begin{equation}\label{directions}
 \g_o^x=\frac{\x_o-\x^\m}{1-\psi^x({\theta}_{\m})}, \quad
 \g_o^y=\frac{\y^\M-\y_o}{\psi^y({\theta}_{\m})-1}.\tag{G2.0}
 \end{equation}
\end{lemma}
\begin{figure}
    \centering
     \scalebox{0.06}
        {\includegraphics{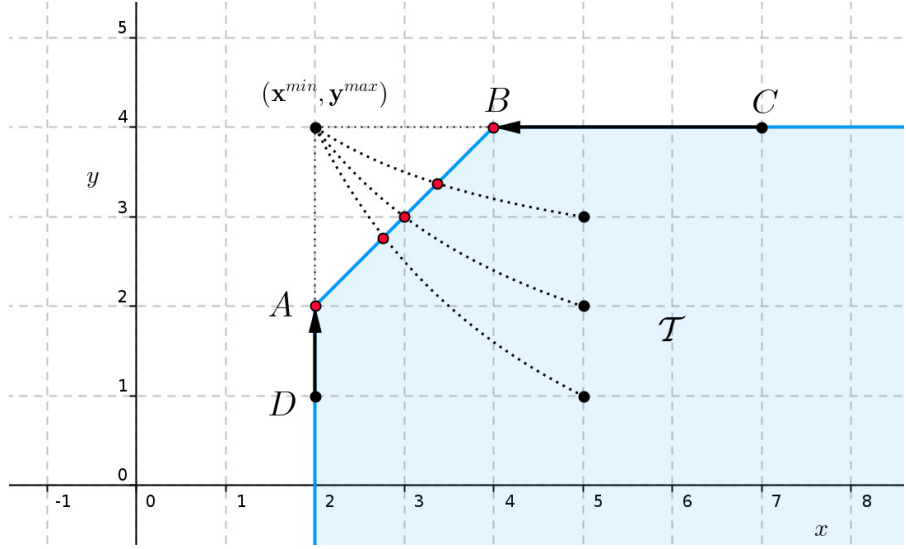}}
    \caption{All  paths with directions \eqref{directions}  intersect at the point $(\x^\m,\y^\M)$ (Lemma~\ref{dir}). Furthermore, \eqref{directions} projects all units in two-dimensional technology $\T$ onto the strongly efficient frontier.}
    \label{F:4}
\end{figure}

\begin{pozn}\label{modif}  The GS models with directions \eqref{directions} satisfy \ref{P5} and \ref{P6} for arbitrary $\psi^x$ and $\psi^y$ with properties \ref{A1}--\ref{A3} by Theorems \ref{invariance} and \ref{translation}. By Lemma \ref{dir}, the same holds for the boundedness property \ref{P4}. Finally, after substituting directions \eqref{directions} into the right-hand sides in \eqref{pathdef} and performing simple calculations, the monotonicity of $\p(\theta)$ obtains for each $\theta\ge\theta_{\min}$. Hence, the monotonicity property \ref{P7}, too, holds universally with directions \eqref{directions}. 

Regarding property \ref{P9}, the situation is more nuanced for two reasons.
First, the super-efficiency considers all points $(\x_o,\y_o)\notin \T$, of which only a small subset with $\x_o\ge \x^{\m}$ and $\y_o\le \y^{\M}$ yields non-negative directions \eqref{directions}. This leads us to consider the absolute values of the expressions in \eqref{directions} as explained in Subsection~\ref{SS5.1}. 
Second, if the range of $\psi^y$ is bounded from below, then the feasibility is compromised on the output side for units (outside $\T$) with very high outputs. Thus, for $\psi^y$ with a bounded range, the super-efficiency fails. In contrast, if the range of $\psi^y$ is unbounded,  by Theorem \ref{superefficiency}\ref{superefficiency.b}, a GS model satisfies super-efficiency for those directions \eqref{directions} that are positive, i.e., super-efficiency holds for all units outside $\T$ except for a negligible set where some direction components are zero.

With respect to properties \ref{P2}, \ref{P3}, and \ref{P8}, the situation is similar to the case of the DDF-\ref{G2} model. The existence of a point $(\x_o,\y_o)\in\partial^W\T$  with $\theta_o^*=1$ will cause all three properties to fail simultaneously (for \ref{P2} and \ref{P3} this is immediate; for \ref{P8} see Theorem \ref{strongmonotonicity}).  Figure \ref{F:3} shows that such points can also be found for directions \eqref{directions} in technological sets of dimension greater than 2. In the case of two-dimensional $\T$, the properties can be satisfied as illustrated in Figure \ref{F:2}  for DDF-\ref{G2} and in Figure \ref{F:4} for linear $\psi^x=\theta$ and hyperbolic $\psi^y=\frac{1}{\theta}$.

\end{pozn}

We now examine the (GS) scheme for linear $\psi^x(\theta)=\theta$ and four choices 
of $\psi^y$ that satisfy the assumptions \ref{A1}--\ref{A3}: $\psi_1^y(\theta)=2-\theta$; $\psi_2^y(\theta)=\theta^{-p}$, $p>0$; $\psi_3^y(\theta)=1-\ln \theta$; and $\psi_4^y(\theta)=\e^{1-\theta}$.
 These functions offer a variety of combinations of domains and ranges, as shown in Table \ref{T:excelent}. Observe that $\psi_1^y$ and $\psi_2^y$ with $p=1$, have already appeared in the context of the DDF-g and HDF-g models, respectively.
To evaluate directions \eqref{directions}, one has to specify the value of $\theta_\m\in[0,1)$. Only for $\psi_1^y$ and $\psi_4^y$ can one choose the ideal value $\theta_\m=0$. The remaining two cases do not have zero in their domain, and for those, we select ad hoc positive values of $\theta_\m$ that yield simple expressions for directions. Intuitively, the choice of $\theta_\m$ affects the ``scaling'' of optimal scores for non-efficient units.
 The resulting values are summarised in Table \ref{T:excelent}.
\begin{table}[t!]
\centering
 \begin{tabular}{l c c c c}
  \hline\\[-1.8ex]
  $\psi^y$ & $2-\theta$& $\theta^{-p},p>0$  & $1-\ln \theta$  &$\e^{1-\theta}$\\[0.8ex]
 \hline\\[-1.8ex]
 domain&$(-\infty,\infty)$&$(0,\infty)$   & $(0,\infty)$  & $(-\infty,\infty)$\\
 image&$(-\infty,\infty)$& $(0,\infty)$   &$(-\infty,\infty)$   & $(0,\infty)$\\
 $\theta_{\m}$&$0$& $1/2$  & $\e^{-1}$  & $0$\\
 $\g_o^x$&$|\x_o-\x^{\m}|$&$2|\x_o-\x^{\m}|$   & $\frac{\e}{\e-1}|\x_o-\x^{\m}|$  &$|\x_o-\x^{\m}|$\\[0.8ex]
  $\g_o^y$&$|\y^{\M}-\y_o|$& $\frac{1}{2^p-1}|\y^{\M}-\y_o|$  &$|\y^{\M}-\y_o|$    &$\frac{1}{\e-1}|\y^{\M}-\y_o|$ \\[0.2ex]
 \ref{P1}, \ref{P4}--\ref{P7}&\checkmark& \checkmark   & \checkmark  &\checkmark\\
  \ref{P2}, \ref{P3}, \ref{P8} &\xmark$^*$&\xmark$^*$   & \xmark$^*$  &\xmark$^*$\\
 \ref{P9}&\checkmark$^{*}$&\xmark   & \checkmark$^{*}$  & \xmark \\
 [0.8ex]
 \hline
 \end{tabular}
 \caption{Values and properties of directions~\eqref{directions} taken in absolute value for $\psi^x(\theta)=\theta$ and three choices of of $\psi^y$ and $\theta_\m$, \checkmark$^{*}$ - the property is satisfied for positive directions; \xmark$^*$ -- the property holds only for some very specific $\T$.}
\label{T:excelent}
\end{table}
Observe that the directions \eqref{directions} for $\psi^y_1$ with the choice of $\theta_{\min}=0$ lead to \ref{G2}. On the contrary, $\psi^y_2$ with $p=1$ and $\theta_{\min}=\frac12$
generate directions whose $x$ component is multiplied by a factor of 2 compared to \ref{G2}.

The statuses with respect to \ref{P1}--\ref{P9} follow from Remark \ref{modif}, where they were analysed for general $\psi^x$ and $\psi^y$. Note that the models with $\psi^y_2$ and $\psi^y_4$ fail \ref{P9} due to the boundedness of their ranges. The property \ref{P10} is not included in Table \ref{T:excelent}, since the results summarised there correspond to \eqref{directions} directions for which the models are not homogeneous; see Theorem~\ref{nonhomog}.

\subsection{Numerical illustration with mixed data}\label{SS5.3}
We present a numerical example that examines the directions \ref{G1} and \eqref{directions}  
for the models with $\psi^x(\theta)=\theta$ and the functions $\psi^y(\theta)$ appearing in Table \ref{T:excelent}. We demonstrate the applicability of our approach to a data-set containing negative data.
To evaluate the efficiency score, as well as to obtain the second phase result, we use the CVX modelling system (see \cite{cvx}, \cite{grant.boyd.08}) to solve convex programmes, implemented in Matlab.
Data for 30 Taiwanese electrical machinery producers are taken from Table 5 of
\cite{tone.al.20} and consist of two inputs and two outputs.  
The minimal and maximal values of the data sample can be seen in Table \ref{T:data}.
The results are presented in Tables \ref{T:example} and \ref{T:example1} in Appendix B.

\begin{table}[H]
\centering
\begin{footnotesize}
 \begin{tabular}{ m{1.6cm}|m{3cm}m{3cm}| m{3cm}m{3cm}}
  &\multicolumn{2}{c|}{Inputs} & \multicolumn{2}{c}{Outputs}  \\
  & Cost of sales & R\&D expenses & Net income & Return  \\
Variable   & $x_1$ & $x_2$ & $y_1$ & $y_2$ \\
    \hline
    \hline
min &	418,242	&  1,319 &  -561,965 & -0.2583 \\
max &	16,655,569	&	634,436 & 3,092,358 & 3.53714 \\
\hline
 \end{tabular}
\end{footnotesize}
 \caption{Minimal and maximal values of the data sample used in the numerical example.}
\label{T:data}
\end{table}

As stated in Section 5.1, the DDF-g model with direction \ref{|G1|} may lead to negative efficiency score, in general. However, this happens only if a negative value 
occurs in each of the inputs in the data-sample. It can be seen from Table \ref{T:data} that the inputs are positive and therefore the efficiency score is bounded bellow by zero.

The choice of a particular function $\psi^y$ has only minimal influence on the ranking of the best-placed units, but the effect is more pronounced for very inefficient units, where the change in ranking can reach 10 places out of 30 companies. On the other hand, the choice of $\psi^y$ has little impact on the number of units projected onto the strongly efficient frontier, where the main determinant is the choice of directions. The number of strongly efficient projections fluctuates between 12 and 13 for the \ref{G1} directions, while it is between 20 and 23 units for the \eqref{directions} directions. The choice of the scaled \ref{G2} direction \eqref{directions} not only leads to a larger number of correct efficiency evaluations, but in the case of the hyperbolic and logarithmic model, it also improves the discriminatory power of the model.

\section{Conclusions}\label{S6}
The paper formulates a general scheme for path-based models that includes BCC, DDF, and HDF as special cases. The models belonging to the scheme are analysed in light of ten desirable properties. Only the unique projection property \ref{P1} is universal; the remaining properties must be examined individually and generally depend on the choice of directions, functions $\psi$, and possibly also on the \emph{configuration of data} that defines the technological set~$\T$. 

The article provides mathematical tools (sufficient and/or necessary conditions to verify specific properties) that allow a detailed analysis of individual models within the GS framework. The usefulness of these tools has been demonstrated on the standard path-based models, on their modifications to general (negative) data, and on models with non-linear and non-hyperbolic functions~$\psi$.

The general scheme is flexible enough to allow modification of existing models designed for positive data to arbitrary (negative) data. For this modification, it is important 
to rewrite models such as BCC-I \eqref{BCCI}, BCC-O \eqref{BCCO}, or HDF \eqref{HDF}, in the form of the GS scheme (Table \ref{Tab:values}), where the expressions for the directions are specified. In the presence of negative data,  these expressions are taken as absolute values. This approach is illustrated with a mixed-data numerical example in Section~\ref{SS5.2}.

The article provides a complete analysis of the  \emph{standard} path-based models and their modifications to negative data with respect to the ten properties (Table \ref{Tab:properties}).  
In addition to the unique projection property \ref{P1}, all standard models satisfy monotonicity \ref{P7}. On the other hand, all models fail \ref{P2}, \ref{P3}, and \ref{P8} simultaneously for \emph{some} $\T$.  Overall, models with \ref{G2} directions satisfy or partially satisfy the highest number of properties, followed by models with \ref{G3} directions.

We have obtained interesting results regarding homogeneity in the general path-based scheme. It is known that in the case of CRS all models in the general scheme meet the property of homogeneity, and thus it could appear that homogeneity and lack of identification are somehow linked. However, the results of this article indicate that this is not the case when the type of technology changes. For VRS, the homogeneity property is satisfied only in a narrow subclass of GS models: In the case of oriented models, homogeneous models are equivalent to BCC-I and BCC-O (Theorem \ref{homogoriented}), and in the case of graph models, only models with directions \ref{G1} and certain special nonlinearities are homogeneous (Theorems \ref{homogT} and \ref{nonhomog}).

 Directions have a strong influence on model properties. Models with positive directions fail \ref{P2},  \ref{P3}, and \ref{P8} for every technology set (Theorems~\ref{weakindication1} and \ref{T:positiveP8} and Remark~\ref{poznG2}). This is especially true for positive data and for all considered directions, except \ref{G2}.  For models with \ref{G2} directions, the status of properties \ref{P2}, \ref{P3}, and \ref{P8} is ambiguous insofar as it depends on the data configuration. There are examples where all three properties fail, but also data configurations (in particular, single-input, single-output technological sets) where all three are satisfied.
 
 The lack of universal validity of \ref{P2}, \ref{P3}, and \ref{P8} can be considered a characteristic feature of path-based models. In contrast, \cite{halicka.trnovska.21} show that the slacks-based models do satisfy these properties. This observation provides a further rationale for distinguishing the two categories of models set out in \cite{russell.schworm.18} -- the path-based versus slacks-based models. We remark that the classification is able to resolve paradoxes in the literature, where some models are currently presented with the adjective `directional' (e.g., \citealp{sahoo.al.14}) suggesting a path-based approach,  despite being in reality slacks-based.

The paper opens up several new research directions.
\begin{itemize}[--]
\item The simultaneous failure of \ref{P2}, \ref{P3}, and \ref{P8} raises the question of whether these properties are linked. It is fairly immediate that \ref{P3} implies  \ref{P2}; such implication is valid also outside our scheme. However, the connection between \ref{P3} and \ref{P8} is not known or obvious. We aim to address this problem in future work.

\item  The general scheme (introduced here in an envelopment form) allows the derivation of a dual multiplicative model, which can then be related to the shadow profit inefficiency. This opens the way for obtaining relationships between technical and overall profit inefficiencies similar to those known for DDF (\citealp{chambers.al.98}).

\item To simplify the exposition, we have only considered technologies with variable returns to scale. The VRS assumption is essential to allow for negative data (see, e.g., \citealp{portela.al.04}). On the other hand, many results in the paper can be adapted to technology sets with other types of returns to scale, at the cost of considering only non-negative data. 

\item The numerical results in Section~\ref{SS5.3} indicate that directions \eqref{directions} may be helpful in practical applications in the sense that a higher number of DMUs are projected onto the strongly efficient frontier compared to other directions. It is natural to ask whether this finding holds up for other data sets.
\end{itemize}

\section*{Acknowledgements }

We would like to thank three anonymous reviewers for detailed comments and constructive suggestions that have lead to several improvements in the paper. The research of the first two authors was supported by the APVV-20-0311 project of the Slovak Research and Development Agency and the VEGA 1/0611/21 grant administered jointly by the Scientific Grant Agency of the Ministry of Education, Science, Research and Sport of the Slovak Republic and the Slovak Academy of Sciences.

\appendix

\section{Proofs of selected theorems}\label{A}

\begin{proof}[Proof of Theorem~\ref{dominance}] 
\ref{T.a} This follows from the assumptions placed on $\g_o$ in \eqref{eq:spc}, the monotonicity properties of $\psi^x$,$\psi^y$ and the definition of $\p_o$ in \eqref{pathdef}.\\[0.2ex]
\ref{T.b}  This follows from \eqref{pathdef} and \ref{A2}.
\\
[0.2ex]\ref{T.c} This part follows from \eqref{pathdef}, \ref{A1}, and \ref{A3}.
\\
[0.2ex]
\ref{T.d} The assumptions on $\g_o$  and part \ref{T.c} of this theorem yield $\lim_{\theta\to a^x_+}\phi_{io}^x(\theta)=-\infty$ for some $i
  $ or $\lim_{\theta\to a^y_+}\phi_{ro}^y(\theta)= +\infty$ for some $r$ .
 From this, part \ref{T.b} of this theorem, and the continuity of $\p_o$, it
   follows that there exists a finite $\bar\theta$ such that $\phi_{ro}^y(\bar\theta)> y^\M_r$ or $\phi_{io}^x(\bar\theta)< x^\m_i$ (see Figure \ref{F:1}). Therefore, by Lemma \ref{L:weak}\ref{W.o} one obtains
    $\p_o(\bar\theta)\notin \T $.
   The property \ref{T.a} of this theorem then shows that for each $\bar\theta$ such that $\p_o(\bar\theta)\notin \T $ one has $\p_o(\theta)\notin \T $ for all $\theta\le \bar\theta$.

On the other hand, $\phi_o(1)=(x_o,y_o)\in \T$ and the property \ref{T.a} of this theorem implies that if $\phi_o(\hat\theta)\in \T $
for some $\hat\theta$, then $\phi_o(\theta)\in \T $ for all $\theta\ge \hat\theta$.
 Now, by the continuity of $\p_o$, there exists a finite $ \theta^*_o$ such that $\theta^*_o=\sup\{\theta : \p_o(\theta)\notin \T\}
 =\inf\{\theta : \p_o(\theta)\in \T\} $. We show that $\p_o(\theta_o^*)$ belongs to the boundary $\partial\T$ of the closed set $\T$.
Assume, in contradiction, that $\p_o(\theta_o^*)\in \text{int}\T$.
From the continuity of $\p_o$  it follows that $\p_o(\theta)\in \text{int}\T$
for all $\theta$ close enough to $\theta_o^*$, and hence also for some $\theta< \theta_o^*$ that contradicts $\theta^*_o =\inf\{\theta : \p_o(\theta)\in \T\}$.
\end{proof}

\begin{proof}[Proof of Theorem~\ref{indication}]
Assume by contradiction that $(\x_o,\y_o)$   is strongly efficient and that there exists an optimal solution $(\la^*_o,\theta^*_o)$  of $(GS)_o$ such that at least one component of the corresponding optimal slacks is positive.  Then the unit $(\X\la^*_o,\Y\la^*_o)$ 
differs from $ (\hat\x_o, \hat\y_o) $, belongs to $\T$  and dominates $(\hat\x_o,\hat\y_o)$, and therefore $(\hat\x_o,\hat\y_o)$ is not strongly efficient.

  To prove the reverse implication,
  we assume that $ (\hat\x_o, \hat\y_o)\in\partial^W\T $. According to Lemma \ref{L:weak}\ref{W.c} there exists $(\bm d^x,\bm d^y)\gneqq\0$,  $(\bm d^x,\bm d^y)\ngtr\0$ such that $(\hat \x_{o}-\bm d^x, \hat\y_o+\bm d^y)\in \T$. By the definition of $\T$,  there then exists $\hat\la\ge\0$  such that $\ones^T\hat\la=1$, $\X\hat\la\le \hat\x_o -\bm d^x$, $\Y\hat\la\ge \hat\y_o+\bm d^y$.
  This implies that $(\hat\la,\theta^*_o)$ is an optimal solution to (GS)$_o$ for $ (\x_o, \y_o) $ and at least one component of inequality \eqref{general2} or \eqref{general3} is strict.
 \end{proof}

\begin{proof}[Proof of Theorem~\ref{invariance}]
The input inequality \eqref{general2} for transformed inputs reads
$$
\bm C\sum_{j=1}^n \x_j\lambda_j=\sum_{j=1}^n(\bm C\x_j)\lambda_j\le \bm C\x_o+(\psi^x(\theta)-1)\g^{Cx}_o \Leftrightarrow \bm C\sum_{j=1}^n\x_j\lambda_j\le \bm C(\x_o+(\psi^y(\theta)-1)\g^{x}_o).$$
Since $\bm C$ is diagonal with positive entries, this is equivalent to the input inequality \eqref{general2} for the original inputs.
This and a similar equivalence for outputs imply that conditions \eqref{general2}--\eqref{general4} before the data transformation are equivalent to those after transformation, and hence the optimal score provided by the model is unchanged.
\end{proof}

\begin{proof} [Proof of Theorem~\ref{translation}] Let $\bm c\in \R^m$ and $\bm b\in \R^s$ are the translation vectors for inputs and outputs, and hence the translated inputs and outputs are   $\x_j+\bm c$ and $\y_j+\bm b$, for all $j=1,\dots,n$ respectively.
The condition $\sum_{j=1}^n\lambda_j=1$ and the invariance of $\g_o^x$ on the translation of inputs imply that the input inequality \eqref{general2} of (GS)$_o$ for translated inputs is equivalent to the input inequality \eqref{general2} of (GS)$_o$ for original inputs:
$$\sum_{j=1}^n(\x_j+\bm c)\lambda_j=\sum_{j=1}^n\x_j\lambda_j+\bm c\le \x_o+\bm c+(\psi^x(\theta)-1)\g_o^x\ \Leftrightarrow \ \sum_{j=1}^n\x_j\lambda_j\le \x_o+(\psi^x(\theta)-1)\g_o^x.$$
Similar equivalence in outputs proves the theorem. \end{proof}

\begin{proof}[Proof of Theorem \ref{superefficiency}] 
\ref{superefficiency.a}  The proof runs along the lines of the proof of Theorem~\ref{dominance}\ref{T.d}. The fact that $\theta^*_o> 1$
    is a simple consequence of the strict monotonicity of $\p_o^x(\theta)$, $\p_o^y(\theta)$, $(\x_o,\y_o)\notin \T$, and $\T$ being closed.\\[0.2ex]
 \ref{superefficiency.b}
 Since $\X\la \le \x^\M$ and $\Y\la\ge \y^\m$,  a sufficient condition for the feasibility of (GS)$_o$ is the existence of $\bar{\theta}\in \D$, such that
\begin{equation}\label{feas1a}
     \x^\M-\x_{o}+\g_{o}^x\le \psi^x(\bar{\theta}) \g_{o}^x,
    \end{equation}
\begin{equation}\label{feas1b}
   \y^\m-\y_o+\g_o^y\ge \psi^y(\bar{\theta}) \g_o^y.
    \end{equation}
Since $\g_o^x>0,$ $\psi^x$ is increasing and
$\lim_{\theta\to+\infty}\psi^x(\theta)=+\infty$, the constraint \eqref{feas1a} is satisfied
for $\bar{\theta}$ sufficiently large.
Since $\psi^y$ is decreasing, the case
$\lim_{\theta\to+\infty}\psi^x(\theta)=+\infty$ implies that the constraint \eqref{feas1b} is satisfied
for $\bar{\theta}$ sufficiently large.  On the other hand,
if $\lim_{\theta\to\infty}\psi^x(\theta)=0$, then the assumptions $\y^\m-\y_o+\g_o^y>0$ and $\g_o^y>0$ guarantee that \eqref{feas1b} is satisfied
for $\bar{\theta}$ sufficiently large. \\[0.2ex]
 \ref{superefficiency.c} Let us first consider the case (i). By assumption, for each feasible $\la$ ($\la\geq0,\ones^T\la=1$), there exists $r$ such that $(\Y\la)_r<y_{ro}-g_{ro}^y$. The assumption on $\psi^y$ implies  $(\Y\la)_r<y_{ro}-g_{ro}^y<y_{ro}-g_{ro}^y+\psi^y(\theta)\g_o^y$ for all $\theta\in\dom(\psi^y)$. Hence we have shown that  $\Y\la\ge \y_o -\g_o^y +\psi^y(\theta)\g_o^y$ is violated  for all  feasible $\la$ and $\theta\in (1,\infty)$ and therefore (GS)$_o$  is infeasible. Consider now the case (ii). Since  $(\Y\la)_r\le \y^{\M}_r$ for all feasible $\la$, we have  $(\Y\la)_r< \y_{ro}-\g_{ro}^y$ and the claim follows by part (i).
\end{proof}
\begin{proof}[Proof of Theorem~\ref{homogT}]
First assume that GS model with \ref{G1} directions is homogeneous of degree $(\alpha,\beta,1)$, i.e. for any technology set 
$\T_+$ the corresponding efficiency measure satisfies $\theta^*(\mu^{\alpha}\x, \mu^{\beta}\y)=\mu \theta^*(\x,\y)$ for each $(x,y)\in\T_+$ and each $\mu>0$. Let $(\x_o,\y_o)\in \mathbb{R}^m_+\times \mathbb{R}^s_+$ and let $\T_+$
be arbitrary, such that $(\x_o,\y_o)\in \partial^S \T_+$. Then any efficiency measure defined over $\T_+$
satisfies $\theta^*(\mu^{\alpha}\x_o, \mu^{\beta}\y_o)=\mu$ for all $\mu>0$ and the projection point for $(\mu^{\alpha}\x_o, \mu^{\beta}\y_o)$ is given by
$$
\p^x(\theta^*(\mu^{\alpha}\x_o, \mu^{\beta}\y_o))=\p^x(\mu)=\psi^x(\mu)\mu^{\alpha}\x_o, \quad
\p^y(\theta^*(\mu^{\alpha}\x_o, \mu^{\beta}\y_o))=\p^y(\mu)=\psi^y(\mu)\mu^{\beta}\y_o.
$$
The set $\mathcal{P}:=\{ (\p^x(\mu), \p^y(\mu)), \ \mu>0\}$ is a smooth path satisfying 
$(\x_o,\y_o)\in \mathcal{P}\subset \partial \T_+$. Recall that $\T_+$ was an arbitrary technology set satisfying $(\x_o,\y_o)\in \partial^S \T_+$. Therefore 
it must hold $ (\p^x(\mu), \p^y(\mu))=(\x_o,\y_o), \forall \mu>0$, which is equivalent to 
$$
\psi^x(\mu)\mu^{\alpha}=1, \quad \psi^y(\mu)\mu^{\beta}=1, \quad \forall \mu>0.
$$
Therefore $\psi^x(\theta)=\theta^{-\alpha}, \psi^y(\theta)=\theta^{-\beta}$;
the ranges for $\alpha, \beta$ follow from \ref{A1bis}.

Now assume that the functions $\psi^x$ and $\psi^y$ are of the form \eqref{homogpath}.
It is easy to see that these functions satisfy both the assumptions \ref{A1}--\ref{A3} and the homogeneity properties: 
\begin{equation}\label{hompsi}
 \psi^x(\theta/\mu)=\mu^{\alpha}\psi^x(\theta), \quad \psi^y(\theta/\mu)=\mu^{\beta}\psi^y(\theta).   
\end{equation}
From \eqref{hompsi} it follows that $\psi^x(\mu)\mu^{\alpha}=1, \psi^y(\mu)\mu^{\beta}=1$ and therefore the efficiency measure 
$$
\theta^*(\mu^{\alpha}\x, \mu^{\beta}\y)=\min \{ \theta\ | \  (\psi^x(\theta)\mu^{\alpha}\x, \psi^y(\theta)\mu^{\beta}\y)\in \T \}
$$
is well-defined provided $(\x,\y)\in \T_+$, since the corresponding minimisation problem is feasible and the optimal value is attained (Theorem \ref{dominance}).
Using \eqref{hompsi} and substitution $\tilde{\theta}:=\theta/\mu$  we get 
\begin{gather*}
\min \{ \theta\ | \  (\psi^x(\theta)\mu^{\alpha}\x, \psi^y(\theta)\mu^{\beta}\y)\in \T \}=\mu\min \{ {\theta}/\mu\ | \  (\psi^x({\theta/\mu})\x, \psi^y({\theta}/\mu)\y)\in \T \}\\
=\mu\min \{ \tilde{\theta}\ | \  (\psi^x(\tilde{\theta})\x, \psi^y(\tilde{\theta})\y)\in \T \}=\mu\theta^*(\x,\y)
\end{gather*}
which proves the homogeneity of $\theta^*$.
\end{proof}

\begin{proof}[Proof of Theorem~\ref{nonhomog}] 
Consider a  GS model with directions $(\g^x,\g^y)=(\mathbf{h}^x,\mathbf{h}^y)$, where 
$(\mathbf{h}^x,\mathbf{h}^y)$ is a constant vector. Assume by contradiction, that the model is homogeneous of degree $(\alpha,\beta,1)$. Similarly, as in the proof of Theorem \ref{homogT} it can be shown that the path
$\mathcal{P}=\{ (\p^x(\mu),\p^y(\mu)), \ \mu>0\}$ defined by
$$
\p^x(\mu)=\mu^{\alpha}\x_o+(\psi^x(\mu)-1)\mathbf{h}^x, \quad
\p^y(\mu)=\mu^{\beta}\y_o+(\psi^y(\mu)-1)\mathbf{h}^y,
$$
satisfies $(\x_o,\y_o)\in \mathcal{P}\subset \partial \T_+$ for any $\T_+$, such that $(\x_o,\y_o)\in \partial^S \T_+$ and hence $ (\p^x(\mu), \p^y(\mu))=(\x_o,\y_o), \forall \mu>0$.
Let $\bar{\mu}>0$ be arbitrary and fixed. Then 
$$
(\bar{\mu}^{\alpha}-1)\x_o=(1-\psi^x(\bar{\mu}))\mathbf{h}^x, \quad
(\bar{\mu}^{\beta}-1)\y_o=(1-\psi^y(\bar{\mu}))\mathbf{h}^y.
$$
While on the right-hand side the vectors are constant, the left-hand side may vary for different $(\x_o,\y_o)$,
which is a contradiction.

The proof for $(\g^x,\g^y)=(\gamma \x_o-\mathbf{h}^x,\mathbf{h}^y-\delta \y_o)$ is analogous.

\end{proof}
\newpage
\section{Tables for the numerical illustrations in Section \ref{SS5.3}}\label{B}
\begin{table}[H]
\centering
\begin{footnotesize}
 \begin{tabular}{ m{1.6cm}|m{1.0cm}m{1.0cm}m{1.2cm}m{1.2cm}| m{1.0cm}m{1.0cm}m{1.2cm}m{1.2cm} }

 Function &\multicolumn{4}{c|}{$\psi^y(\theta)=2-\theta$} & \multicolumn{4}{c}{$\psi^y(\theta)=1/\theta$}  \\
 Directions & \multicolumn{2}{c}{\ref{G1}} & \multicolumn{2}{l|}{\eqref{directions} $\theta_{\m}=0$} &\multicolumn{2}{c}{\ref{G1}} &\multicolumn{2}{l}{\eqref{directions} $\theta_{\m}=0.01$ } \\
 DMU & Score & Rank & Score & Rank & Score & Rank & Score & Rank \\
    \hline
    \hline
1	&	0.66*	&	15	&	0.759	&	28	&	0.713*	&	14	&	0.401	&	14	\\
2	&	1*	&	1-9	&	1*	&	1-9	&	1*	&	1-9	&	1*	&	1-9	\\
3	&	0.683	&	14	&	0.888	&	12	&	0.691	&	15	&	0.578*	&	11	\\
4	&	1*	&	1-9	&	1*	&	1-9	&	1*	&	1-9	&	1*	&	1-9	\\
5	&	1*	&	1-9	&	1*	&	1-9	&	1*	&	1-9	&	1*	&	1-9	\\
6	&	0.175	&	29	&	0.784*	&	26	&	0.256	&	29	&	0.145*	&	28	\\
7	&	0.69	&	13	&	0.86*	&	17	&	0.763	&	13	&	0.254*	&	19	\\
8	&	0.283	&	26	&	0.801	&	24	&	0.3	&	28	&	0.167*	&	24	\\
9	&	0.126*	&	30	&	0.744*	&	29	&	0.194*	&	30	&	0.112*	&	30	\\
10	&	0.372	&	22	&	0.777	&	27	&	0.526	&	19	&	0.156*	&	27	\\
11	&	0.37	&	23	&	0.707*	&	30	&	0.512	&	20	&	0.169*	&	23	\\
12	&	1*	&	1-9	&	1*	&	1-9	&	1*	&	1-9	&	1*	&	1-9	\\
13	&	0.927	&	10	&	0.968*	&	10	&	0.929	&	10	&	0.876	&	10	\\
14	&	1*	&	1-9	&	1*	&	1-9	&	1*	&	1-9	&	1*	&	1-9	\\
15	&	0.217	&	28	&	0.793*	&	25	&	0.339	&	25	&	0.114	&	29	\\
16	&	0.658	&	16	&	0.881	&	15	&	0.68	&	16	&	0.309	&	16	\\
17	&	0.22	&	27	&	0.807*	&	22	&	0.319	&	27	&	0.165*	&	25	\\
18	&	0.438	&	19	&	0.861	&	16	&	0.494	&	21	&	0.233*	&	21	\\
19	&	1*	&	1-9	&	1*	&	1-9	&	1*	&	1-9	&	1*	&	1-9	\\
20	&	0.415	&	21	&	0.841	&	18	&	0.448	&	22	&	0.165	&	26	\\
21	&	0.335	&	24	&	0.804	&	23	&	0.408	&	24	&	0.27	&	17	\\
22	&	1*	&	1-9	&	1*	&	1-9	&	1*	&	1-9	&	1*	&	1-9	\\
23	&	1*	&	1-9	&	1*	&	1-9	&	1*	&	1-9	&	1*	&	1-9	\\
24	&	0.765	&	12	&	0.882	&	14	&	0.765	&	12	&	0.547	&	13	\\
25	&	0.431	&	20	&	0.814*	&	20	&	0.442	&	23	&	0.372	&	15	\\
26	&	1*	&	1-9	&	1*	&	1-9	&	1*	&	1-9	&	1*	&	1-9	\\
27	&	0.305	&	25	&	0.815*	&	19	&	0.33	&	26	&	0.255	&	18	\\
28	&	0.447	&	18	&	0.812*	&	21	&	0.581	&	18	&	0.176*	&	22	\\
29	&	0.835*	&	11	&	0.933*	&	11	&	0.85*	&	11	&	0.576*	&	12	\\
30	&	0.452	&	17	&	0.883*	&	13	&	0.598*	&	17	&	0.239*	&	20	\\
\hline
Average & 0.627 & & 0.880 & & 0.671 & & 0.509 & \\
Minimum & 0.126 & & 0.707 & & 0.194 & & 0.112 & \\
Correct & 12 & & 21 & & 13 & & 21 & \\
 \end{tabular}
\end{footnotesize}
 \caption{Efficiency scores for models with $\psi^x=\theta$, two choices of $\psi^y$, directions \ref{G1} and  \eqref{directions}.  An asterisk next to the efficiency value indicates projection onto the strongly efficient frontier. The last row indicates number of units projected onto the strongly efficient frontier.}
\label{T:example}
\end{table}
\begin{table}[H]
\centering
\begin{footnotesize}
 \begin{tabular}{ m{1.6cm}|m{1cm}m{1.0cm}m{1.2cm}m{1.2cm}| m{1cm}m{1.0cm}m{1.2cm}m{1.2cm} }

 Function &\multicolumn{4}{c|}{$\psi^y(\theta)=\exp(1-\theta)$} & \multicolumn{4}{c}{$\psi^y(\theta)=1-\ln\theta$}  \\
 Directions & \multicolumn{2}{c}{\ref{G1}} & \multicolumn{2}{l|}{\eqref{directions} $\theta_{\m}=0$} &\multicolumn{2}{c}{\ref{G1}} &\multicolumn{2}{l}{\eqref{directions} $\theta_{\m}=0.01$ } \\
 DMU & Score & Rank & Score & Rank & Score & Rank & Score & Rank \\
    \hline
    \hline
1	&	0,686*	&	14	&	0,697	&	28	& 0,69* &	14	&	0,565*	&	18	\\
2	&	1*	&	1-9	&	1*	&	1	&	1*	&	1	&	1*	&	1	\\
3	&	0,686	&	15	&	0,833	&	12	&	0,687	&	15	&	0,686*	&	12	\\
4	&	1*	&	1-9	&	1*	&	2	&	1*	&	2	&	1*	&	2	\\
5	&	1*	&	1-9	&	1*	&	3	&	1*	&	3	&	1*	&	3	\\
6	&	0,19	&	29	&	0,698*	&	26	&	0,203	&	29	&	0,454*	&	27	\\
7	&	0,73	&	13	&	0,799*	&	16	&	0,733	&	13	&	0,627*	&	15	\\
8	&	0,283	&	26	&	0,711	&	23	&	0,283	&	26	&	0,437	&	28	\\
9	&	0,131*	&	30	&	0,652*	&	29	&	0,137*	&	30	&	0,407*	&	30	\\
10	&	0,439	&	20	&	0,698	&	27	&	0,457	&	20	&	0,5*	&	23	\\
11	&	0,429	&	22	&	0,626*	&	30	&	0,446	&	21	&	0,422	&	29	\\
12	&	1*	&	1-9	&	1*	&	4	&	1*	&	4	&	1*	&	4	\\
13	&	0,928	&	10	&	0,953	&	10	&	0,928	&	10	&	0,922	&	10	\\
14	&	1*	&	1-9	&	1*	&	5	&	1*	&	5	&	1*	&	5	\\
15	&	0,232	&	28	&	0,708*	&	25	&	0,246	&	28	&	0,464*	&	26	\\
16	&	0,668	&	16	&	0,827	&	13	&	0,669	&	16	&	0,678	&	13	\\
17	&	0,238	&	27	&	0,726*	&	22	&	0,251	&	27	&	0,488*	&	24	\\
18	&	0,457	&	19	&	0,79	&	17	&	0,463	&	19	&	0,566	&	17	\\
19	&	1*	&	1-9	&	1*	&	6	&	1*	&	6	&	1*	&	6	\\
20	&	0,425	&	23	&	0,763	&	18	&	0,428	&	23	&	0,512	&	22	\\
21	&	0,357	&	24	&	0,711	&	24	&	0,366	&	24	&	0,527*	&	21	\\
22	&	1*	&	1-9	&	1*	&	7	&	1*	&	7	&	1*	&	7	\\
23	&	1*	&	1-9	&	1*	&	8	&	1*	&	8	&	1*	&	8	\\
24	&	0,765	&	12	&	0,82	&	15	&	0,765	&	12	&	0,618	&	16	\\
25	&	0,432	&	21	&	0,729*	&	21	&	0,434	&	22	&	0,543*	&	20	\\
26	&	1*	&	1-9	&	1*	&	9	&	1*	&	9	&	1*	&	9	\\
27	&	0,308	&	25	&	0,733*	&	20	&	0,31	&	25	&	0,484*	&	25	\\
28	&	0,51	&	18	&	0,741*	&	19	&	0,523	&	18	&	0,551*	&	19	\\
29	&	0,842*	&	11	&	0,899*	&	11	&	0,843*	&	11	&	0,808*	&	11	\\
30	&	0,524	&	17	&	0,823*	&	14	&	0,537	&	17	&	0,64*	&	14	\\

\hline
Average & 0.642 & & 0.831 & & 0.647 & & 0.697 & \\
Minimum & 0.131 & & 0.626 & & 0.137 & & 0.407 & \\
Correct & 12 & & 20 & & 12 & & 23 & \\
 \end{tabular}
\end{footnotesize}
 \caption{Efficiency scores for models with $\psi^x=\theta$, two choices of $\psi^y$, directions \ref{G1} and  \eqref{directions}.  An asterisk next to the efficiency value indicates projection onto the strongly efficient frontier. The last row indicates number of units projected onto the strongly efficient frontier.}
\label{T:example1}
\end{table}


\newpage


{\small
\begin{thebibliography}{46}
\expandafter\ifx\csname natexlab\endcsname\relax\def\natexlab#1{#1}\fi
\providecommand{\url}[1]{\texttt{#1}}
\providecommand{\href}[2]{#2}
\providecommand{\path}[1]{#1}
\providecommand{\DOIprefix}{doi:}
\providecommand{\ArXivprefix}{arXiv:}
\providecommand{\URLprefix}{URL: }
\providecommand{\Pubmedprefix}{pmid:}
\providecommand{\doi}[1]{\href{http://dx.doi.org/#1}{\path{#1}}}
\providecommand{\Pubmed}[1]{\href{pmid:#1}{\path{#1}}}
\providecommand{\bibinfo}[2]{#2}
\ifx\xfnm\relax \def\xfnm[#1]{\unskip,\space#1}\fi
\bibitem[{Aparicio and Monge(2022)}]{aparicio.monge.22}
\bibinfo{author}{Aparicio, J.}, \bibinfo{author}{Monge, J.F.},
  \bibinfo{year}{2022}.
\newblock \bibinfo{title}{The generalized range adjusted measure in data
  envelopment analysis: properties, computational aspects and duality}.
\newblock \bibinfo{journal}{European J. Oper. Res.} \bibinfo{volume}{302},
  \bibinfo{pages}{621--632}.
\newblock \DOIprefix\doi{10.1016/j.ejor.2022.01.001}.
\bibitem[{Aparicio et~al.(2013)Aparicio, Pastor and Ray}]{aparicio.al.13}
\bibinfo{author}{Aparicio, J.}, \bibinfo{author}{Pastor, J.T.},
  \bibinfo{author}{Ray, S.C.}, \bibinfo{year}{2013}.
\newblock \bibinfo{title}{An overall measure of technical inefficiency at the
  firm and at the industry level: {T}he `lost profit on outlay'}.
\newblock \bibinfo{journal}{European J. Oper. Res.} \bibinfo{volume}{226},
  \bibinfo{pages}{154--162}.
\newblock \DOIprefix\doi{10.1016/j.ejor.2012.10.028}.
\bibitem[{Aparicio et~al.(2016)Aparicio, Pastor and Vidal}]{aparicio.al.16}
\bibinfo{author}{Aparicio, J.}, \bibinfo{author}{Pastor, J.T.},
  \bibinfo{author}{Vidal, F.}, \bibinfo{year}{2016}.
\newblock \bibinfo{title}{The directional distance function and the translation
  invariance property}.
\newblock \bibinfo{journal}{Omega} \bibinfo{volume}{58}, \bibinfo{pages}{1--3}.
\newblock \DOIprefix\doi{10.1016/j.omega.2015.04.012}.
\bibitem[{Banker et~al.(1984)Banker, Charnes and Cooper}]{banker.al.84}
\bibinfo{author}{Banker, R.D.}, \bibinfo{author}{Charnes, A.},
  \bibinfo{author}{Cooper, W.W.}, \bibinfo{year}{1984}.
\newblock \bibinfo{title}{Some models for estimating technical and scale
  inefficiencies in data envelopment analysis}.
\newblock \bibinfo{journal}{Management Sci.} \bibinfo{volume}{30},
  \bibinfo{pages}{1078--1092}.
\newblock \DOIprefix\doi{10.1287/mnsc.30.9.1078}.
\bibitem[{Briec(1997)}]{briec.97}
\bibinfo{author}{Briec, W.}, \bibinfo{year}{1997}.
\newblock \bibinfo{title}{A graph-type extension of {F}arrell technical
  efficiency measure}.
\newblock \bibinfo{journal}{J. Prod. Anal.} \bibinfo{volume}{8},
  \bibinfo{pages}{95--110}.
\newblock \DOIprefix\doi{10.1023/a:1007728515733}.
\bibitem[{Briec and Kerstens(2009)}]{briec.kerstens.09}
\bibinfo{author}{Briec, W.}, \bibinfo{author}{Kerstens, K.},
  \bibinfo{year}{2009}.
\newblock \bibinfo{title}{Infeasibility and directional distance functions with
  application to the determinateness of the {L}uenberger productivity
  indicator}.
\newblock \bibinfo{journal}{J. Optim. Theory Appl.} \bibinfo{volume}{141},
  \bibinfo{pages}{55--73}.
\newblock \DOIprefix\doi{10.1007/s10957-008-9503-2}.
\bibitem[{Briec and Lesourd(1999)}]{briec.lesourd.99}
\bibinfo{author}{Briec, W.}, \bibinfo{author}{Lesourd, J.B.},
  \bibinfo{year}{1999}.
\newblock \bibinfo{title}{Metric distance function and profit: {S}ome duality
  results}.
\newblock \bibinfo{journal}{J. Optim. Theory Appl.} \bibinfo{volume}{101},
  \bibinfo{pages}{15--33}.
\newblock \DOIprefix\doi{10.1023/A:1021762809393}.
\bibitem[{Chambers et~al.(1996a)Chambers, Chung and
  F{\"a}re}]{chambers.al.96.jet}
\bibinfo{author}{Chambers, R.G.}, \bibinfo{author}{Chung, Y.},
  \bibinfo{author}{F{\"a}re, R.}, \bibinfo{year}{1996}a.
\newblock \bibinfo{title}{Benefit and distance functions}.
\newblock \bibinfo{journal}{J. Econom. Theory} \bibinfo{volume}{70},
  \bibinfo{pages}{407--419}.
\newblock \DOIprefix\doi{10.1006/jeth.1996.0096}.
\bibitem[{Chambers et~al.(1998)Chambers, Chung and F\"{a}re}]{chambers.al.98}
\bibinfo{author}{Chambers, R.G.}, \bibinfo{author}{Chung, Y.},
  \bibinfo{author}{F\"{a}re, R.}, \bibinfo{year}{1998}.
\newblock \bibinfo{title}{Profit, directional distance functions, and
  {N}erlovian efficiency}.
\newblock \bibinfo{journal}{J. Optim. Theory Appl.} \bibinfo{volume}{98},
  \bibinfo{pages}{351--364}.
\newblock \DOIprefix\doi{10.1023/A:1022637501082}.
\bibitem[{Chambers et~al.(1996b)Chambers, F{\=a}ure and
  Grosskopf}]{chambers.al.96.per}
\bibinfo{author}{Chambers, R.G.}, \bibinfo{author}{F{\=a}ure, R.},
  \bibinfo{author}{Grosskopf, S.}, \bibinfo{year}{1996}b.
\newblock \bibinfo{title}{Productivity growth in {APEC} countries}.
\newblock \bibinfo{journal}{Pacific Econ. Rev.} \bibinfo{volume}{1},
  \bibinfo{pages}{181--190}.
\newblock \DOIprefix\doi{10.1111/j.1468-0106.1996.tb00184.x}.
\bibitem[{Charnes et~al.(1985)Charnes, Cooper, Golany, Seiford and
  Stutz}]{charnes.al.85}
\bibinfo{author}{Charnes, A.}, \bibinfo{author}{Cooper, W.W.},
  \bibinfo{author}{Golany, B.}, \bibinfo{author}{Seiford, L.},
  \bibinfo{author}{Stutz, J.}, \bibinfo{year}{1985}.
\newblock \bibinfo{title}{Foundations of data envelopment analysis for
  {P}areto-{K}oopmans efficient empirical production functions}.
\newblock \bibinfo{journal}{J. Econometrics} \bibinfo{volume}{30},
  \bibinfo{pages}{91--107}.
\newblock \DOIprefix\doi{10.1016/0304-4076(85)90133-2}.
\bibitem[{Chavas and Cox(1999)}]{chavas.cox.99}
\bibinfo{author}{Chavas, J.P.}, \bibinfo{author}{Cox, T.L.},
  \bibinfo{year}{1999}.
\newblock \bibinfo{title}{A generalized distance function and the analysis of
  production efficiency}.
\newblock \bibinfo{journal}{South. Econ. J.} \bibinfo{volume}{66},
  \bibinfo{pages}{294--318}.
\newblock \DOIprefix\doi{https://doi.org/10.1002/j.2325-8012.1999.tb00248.x}.
\bibitem[{Cheng et~al.(2013)Cheng, Zervopoulos and Qian}]{cheng.al.13}
\bibinfo{author}{Cheng, G.}, \bibinfo{author}{Zervopoulos, P.},
  \bibinfo{author}{Qian, Z.}, \bibinfo{year}{2013}.
\newblock \bibinfo{title}{A variant of radial measure capable of dealing with
  negative inputs and outputs in data envelopment analysis}.
\newblock \bibinfo{journal}{European J. Oper. Res.} \bibinfo{volume}{225},
  \bibinfo{pages}{100--105}.
\newblock \DOIprefix\doi{10.1016/j.ejor.2012.09.031}.
\bibitem[{Cooper et~al.(1999)Cooper, Park and Pastor}]{cooper.al.99}
\bibinfo{author}{Cooper, W.W.}, \bibinfo{author}{Park, K.S.},
  \bibinfo{author}{Pastor, J.T.}, \bibinfo{year}{1999}.
\newblock \bibinfo{title}{{RAM}: A range adjusted measure of inefficiency for
  use with additive models, and relations to other models and measures in
  {DEA}}.
\newblock \bibinfo{journal}{J. Prod. Anal.} \bibinfo{volume}{11},
  \bibinfo{pages}{5--42}.
\newblock \DOIprefix\doi{10.1023/a:1007701304281}.
\bibitem[{Cooper et~al.(2011)Cooper, Pastor, Borras, Aparicio and
  Pastor}]{cooper.11}
\bibinfo{author}{Cooper, W.W.}, \bibinfo{author}{Pastor, J.T.},
  \bibinfo{author}{Borras, F.}, \bibinfo{author}{Aparicio, J.},
  \bibinfo{author}{Pastor, D.}, \bibinfo{year}{2011}.
\newblock \bibinfo{title}{{BAM}: A bounded adjusted measure of efficiency for
  use with bounded additive models}.
\newblock \bibinfo{journal}{J. Prod. Anal.} \bibinfo{volume}{35},
  \bibinfo{pages}{85--94}.
\newblock \DOIprefix\doi{10.1007/s11123-010-0190-2}.
\bibitem[{Cooper et~al.(2007)Cooper, Seiford and Tone}]{cooper.al.07}
\bibinfo{author}{Cooper, W.W.}, \bibinfo{author}{Seiford, L.M.},
  \bibinfo{author}{Tone, K.}, \bibinfo{year}{2007}.
\newblock \bibinfo{title}{Data Envelopment Analysis: {A} Comprehensive Text
  with Models, Applications, References and Dea-Solver Software}.
\newblock \bibinfo{publisher}{Springer {US}}.
\newblock \DOIprefix\doi{10.1007/978-0-387-45283-8}.
\bibitem[{Cuesta and Zof{\'\i}o(2005)}]{cuesta.zofio.05}
\bibinfo{author}{Cuesta, R.A.}, \bibinfo{author}{Zof{\'\i}o, J.L.},
  \bibinfo{year}{2005}.
\newblock \bibinfo{title}{Hyperbolic efficiency and parametric distance
  functions: {W}ith application to {S}panish savings banks}.
\newblock \bibinfo{journal}{J. Prod. Anal.} \bibinfo{volume}{24},
  \bibinfo{pages}{31--48}.
\newblock \DOIprefix\doi{10.1007/s11123-005-3039-3}.
\bibitem[{F\"{a}re et~al.(1985)F\"{a}re, Grosskopf and Lovell}]{fare.al.85}
\bibinfo{author}{F\"{a}re, R.}, \bibinfo{author}{Grosskopf, S.},
  \bibinfo{author}{Lovell, C.A.K.}, \bibinfo{year}{1985}.
\newblock \bibinfo{title}{The Measurement of Efficiency of Production}.
\newblock \bibinfo{publisher}{Springer Netherlands}.
\newblock \DOIprefix\doi{10.1007/978-94-015-7721-2}.
\bibitem[{F{\"{a}}re et~al.(2008)F{\"{a}}re, Grosskopf and
  Margaritis}]{fare.al.08}
\bibinfo{author}{F{\"{a}}re, R.}, \bibinfo{author}{Grosskopf, S.},
  \bibinfo{author}{Margaritis, D.}, \bibinfo{year}{2008}.
\newblock \bibinfo{title}{Efficiency and productivity: {M}almquist and more},
  in: \bibinfo{editor}{Fried, H.O.}, \bibinfo{editor}{Knox~Lovell, C.A.},
  \bibinfo{editor}{Schmidt, S.S.} (Eds.), \bibinfo{booktitle}{The Measurement
  of Productive Efficiency and Productivity Growth}. \bibinfo{publisher}{Oxford
  University Press}, pp. \bibinfo{pages}{522--622}.
\newblock \DOIprefix\doi{10.1093/acprof:oso/9780195183528.003.0005}.
\bibitem[{F\"{a}re and Lovell(1978)}]{fare.lovell.78}
\bibinfo{author}{F\"{a}re, R.}, \bibinfo{author}{Lovell, C.A.K.},
  \bibinfo{year}{1978}.
\newblock \bibinfo{title}{Measuring the technical efficiency of production}.
\newblock \bibinfo{journal}{J. Econom. Theory} \bibinfo{volume}{19},
  \bibinfo{pages}{150--162}.
\newblock \DOIprefix\doi{10.1016/0022-0531(78)90060-1}.
\bibitem[{F{\"a}re et~al.(1983)F{\"a}re, Lovell and Zieschang}]{fare.al.83}
\bibinfo{author}{F{\"a}re, R.}, \bibinfo{author}{Lovell, C.A.K.},
  \bibinfo{author}{Zieschang, K.}, \bibinfo{year}{1983}.
\newblock \bibinfo{title}{Measuring the technical efficiency of multiple output
  production technologies}, in: \bibinfo{editor}{Eichhorn, W.},
  \bibinfo{editor}{Henn, R.}, \bibinfo{editor}{Neumann, K.},
  \bibinfo{editor}{Shephard, R.W.} (Eds.), \bibinfo{booktitle}{Quantitative
  Studies on Production and Prices}, \bibinfo{publisher}{Physica-Verlag},
  \bibinfo{address}{Heidelberg}. pp. \bibinfo{pages}{159--171}.
\newblock \DOIprefix\doi{10.1007/978-3-662-41526-9_12}.
\bibitem[{F\"{a}re et~al.(2016)F\"{a}re, Margaritis, Rouse and
  Roshdi}]{fare.al.16}
\bibinfo{author}{F\"{a}re, R.}, \bibinfo{author}{Margaritis, D.},
  \bibinfo{author}{Rouse, P.}, \bibinfo{author}{Roshdi, I.},
  \bibinfo{year}{2016}.
\newblock \bibinfo{title}{Estimating the hyperbolic distance function: {A}
  directional distance function approach}.
\newblock \bibinfo{journal}{European J. Oper. Res.} \bibinfo{volume}{254},
  \bibinfo{pages}{312--319}.
\newblock \DOIprefix\doi{10.1016/j.ejor.2016.03.045}.
\bibitem[{Farrell(1957)}]{farrell.57}
\bibinfo{author}{Farrell, M.J.}, \bibinfo{year}{1957}.
\newblock \bibinfo{title}{The measurement of productive efficiency}.
\newblock \bibinfo{journal}{J. Roy. Statist. Soc. Ser. A}
  \bibinfo{volume}{120}, \bibinfo{pages}{253--281}.
\newblock \DOIprefix\doi{10.2307/2343100}.
\bibitem[{Grant and Boyd(2008)}]{grant.boyd.08}
\bibinfo{author}{Grant, M.}, \bibinfo{author}{Boyd, S.}, \bibinfo{year}{2008}.
\newblock \bibinfo{title}{Graph implementations for nonsmooth convex programs},
  in: \bibinfo{editor}{Blondel, V.}, \bibinfo{editor}{Boyd, S.},
  \bibinfo{editor}{Kimura, H.} (Eds.), \bibinfo{booktitle}{Recent Advances in
  Learning and Control}. \bibinfo{publisher}{Springer-Verlag Limited}. Lecture
  Notes in Control and Information Sciences, pp. \bibinfo{pages}{95--110}.
\newblock \DOIprefix\doi{10.1007/978-1-84800-155-8_7}.
\bibitem[{Grant and Boyd(2014)}]{cvx}
\bibinfo{author}{Grant, M.}, \bibinfo{author}{Boyd, S.}, \bibinfo{year}{2014}.
\newblock \bibinfo{title}{{CVX}: Matlab software for disciplined convex
  programming, version 2.1}.
\newblock \bibinfo{howpublished}{\url{http://cvxr.com/cvx}}.
\bibitem[{Halick\'{a} and Trnovsk\'{a}(2019)}]{halicka.trnovska.19}
\bibinfo{author}{Halick\'{a}, M.}, \bibinfo{author}{Trnovsk\'{a}, M.},
  \bibinfo{year}{2019}.
\newblock \bibinfo{title}{Duality and profit efficiency for the hyperbolic
  measure model}.
\newblock \bibinfo{journal}{European J. Oper. Res.} \bibinfo{volume}{278},
  \bibinfo{pages}{410--421}.
\newblock \DOIprefix\doi{10.1016/j.ejor.2018.12.001}.
\bibitem[{Halick\'{a} and Trnovsk\'{a}(2021)}]{halicka.trnovska.21}
\bibinfo{author}{Halick\'{a}, M.}, \bibinfo{author}{Trnovsk\'{a}, M.},
  \bibinfo{year}{2021}.
\newblock \bibinfo{title}{A unified approach to non-radial graph models in data
  envelopment analysis: {C}ommon features, geometry, and duality}.
\newblock \bibinfo{journal}{European J. Oper. Res.} \bibinfo{volume}{289},
  \bibinfo{pages}{611--627}.
\newblock \DOIprefix\doi{10.1016/j.ejor.2020.07.019}.
\bibitem[{Hasannasab et~al.(2019)Hasannasab, Margaritis, Roshdi and
  Rouse}]{hasannasab.al.19}
\bibinfo{author}{Hasannasab, M.}, \bibinfo{author}{Margaritis, D.},
  \bibinfo{author}{Roshdi, I.}, \bibinfo{author}{Rouse, P.},
  \bibinfo{year}{2019}.
\newblock \bibinfo{title}{Hyperbolic efficiency measurement: {A} conic
  programming approach}.
\newblock \bibinfo{journal}{European J. Oper. Res.} \bibinfo{volume}{278},
  \bibinfo{pages}{401--409}.
\newblock \DOIprefix\doi{10.1016/j.ejor.2018.12.005}.
\bibitem[{Hudgins and Primont(2007)}]{hudgins.primont.07}
\bibinfo{author}{Hudgins, L.B.}, \bibinfo{author}{Primont, D.},
  \bibinfo{year}{2007}.
\newblock \bibinfo{title}{Derivative properties of directional technology
  distance functions}, in: \bibinfo{editor}{F{\"a}re, R.},
  \bibinfo{editor}{Grosskopf, S.}, \bibinfo{editor}{Primont, D.} (Eds.),
  \bibinfo{booktitle}{Aggregation, Efficiency, and Measurement}.
  \bibinfo{publisher}{Springer, Boston, MA}, pp. \bibinfo{pages}{31--43}.
\newblock \DOIprefix\doi{10.1007/978-0-387-47677-3_2}.
\bibitem[{Johnson and McGinnis(2009)}]{johnson.mcginnis.09}
\bibinfo{author}{Johnson, A.L.}, \bibinfo{author}{McGinnis, L.F.},
  \bibinfo{year}{2009}.
\newblock \bibinfo{title}{The hyperbolic-oriented efficiency measure as a
  remedy to infeasibility of super efficiency models}.
\newblock \bibinfo{journal}{J. Oper. Res. Soc.} \bibinfo{volume}{60},
  \bibinfo{pages}{1511--1517}.
\newblock \DOIprefix\doi{10.1057/jors.2009.71}.
\bibitem[{Kerstens and Van~de Woestyne(2011)}]{kerstens.vandewoestyne.11}
\bibinfo{author}{Kerstens, K.}, \bibinfo{author}{Van~de Woestyne, I.},
  \bibinfo{year}{2011}.
\newblock \bibinfo{title}{Negative data in {DEA}: {A} simple proportional
  distance function approach}.
\newblock \bibinfo{journal}{J. Oper. Res. Soc.} \bibinfo{volume}{62},
  \bibinfo{pages}{1413--1419}.
\newblock \DOIprefix\doi{10.1057/jors.2010.108}.
\bibitem[{Levkoff et~al.(2011)Levkoff, Russell and Schworm}]{levkoff.al.11}
\bibinfo{author}{Levkoff, S.B.}, \bibinfo{author}{Russell, R.R.},
  \bibinfo{author}{Schworm, W.}, \bibinfo{year}{2011}.
\newblock \bibinfo{title}{Boundary problems with the
  {\textquotedblleft}{R}ussell{\textquotedblright} graph measure of technical
  efficiency: {A} refinement}.
\newblock \bibinfo{journal}{J. Prod. Anal.} \bibinfo{volume}{37},
  \bibinfo{pages}{239--248}.
\newblock \DOIprefix\doi{10.1007/s11123-011-0241-3}.
\bibitem[{Mehdiloozad and Roshdi(2014)}]{mehdiloozad.roshdi.14.arxiv}
\bibinfo{author}{Mehdiloozad, M.}, \bibinfo{author}{Roshdi, I.},
  \bibinfo{year}{2014}.
\newblock \bibinfo{title}{Analyzing the concept of super-efficiency in data
  envelopment analysis: {A} directional distance function approach}.
\newblock \DOIprefix\doi{10.48550/arXiv.1407.2599}.
\bibitem[{Pastor et~al.(2022)Pastor, Aparicio and Zof{\'\i}o}]{pastor.al.22}
\bibinfo{author}{Pastor, J.T.}, \bibinfo{author}{Aparicio, J.},
  \bibinfo{author}{Zof{\'\i}o, J.L.}, \bibinfo{year}{2022}.
\newblock \bibinfo{title}{Benchmarking Economic Efficiency: Technical and
  Allocative Fundamentals}.
\newblock International Series in Operations Research and Management Science,
  \bibinfo{publisher}{Springer Cham}.
\newblock \DOIprefix\doi{10.1007/978-3-030-84397-7}.
\bibitem[{Pastor et~al.(1999)Pastor, Ruiz and Sirvent}]{pastor.al.99}
\bibinfo{author}{Pastor, J.T.}, \bibinfo{author}{Ruiz, J.L.},
  \bibinfo{author}{Sirvent, I.}, \bibinfo{year}{1999}.
\newblock \bibinfo{title}{An enhanced {DEA} {R}ussell graph efficiency
  measure}.
\newblock \bibinfo{journal}{European J. Oper. Res.} \bibinfo{volume}{115},
  \bibinfo{pages}{596--607}.
\newblock \DOIprefix\doi{10.1016/S0377-2217(98)00098-8}.
\bibitem[{Portela et~al.(2004)Portela, Thanassoulis and
  Simpson}]{portela.al.04}
\bibinfo{author}{Portela, M.C.A.S.}, \bibinfo{author}{Thanassoulis, E.},
  \bibinfo{author}{Simpson, G.}, \bibinfo{year}{2004}.
\newblock \bibinfo{title}{Negative data in {DEA}: {A} directional distance
  approach applied to bank branches}.
\newblock \bibinfo{journal}{J. Oper. Res. Soc.} \bibinfo{volume}{55},
  \bibinfo{pages}{1111--1121}.
\newblock \DOIprefix\doi{10.1057/palgrave.jors.2601768}.
\bibitem[{Roshdi et~al.(2018)Roshdi, Hasannasab, Margaritis and
  Rouse}]{roshdi.al.18}
\bibinfo{author}{Roshdi, I.}, \bibinfo{author}{Hasannasab, M.},
  \bibinfo{author}{Margaritis, D.}, \bibinfo{author}{Rouse, P.},
  \bibinfo{year}{2018}.
\newblock \bibinfo{title}{Generalised weak disposability and efficiency
  measurement in environmental technologies}.
\newblock \bibinfo{journal}{European J. Oper. Res.} \bibinfo{volume}{266},
  \bibinfo{pages}{1000--1012}.
\newblock \DOIprefix\doi{10.1016/j.ejor.2017.10.033}.
\bibitem[{Russell and Schworm(2008)}]{russell.schworm.08}
\bibinfo{author}{Russell, R.R.}, \bibinfo{author}{Schworm, W.},
  \bibinfo{year}{2008}.
\newblock \bibinfo{title}{Axiomatic foundations of efficiency measurement on
  data-generated technologies}.
\newblock \bibinfo{journal}{J. Prod. Anal.} \bibinfo{volume}{31},
  \bibinfo{pages}{77--86}.
\newblock \DOIprefix\doi{10.1007/s11123-008-0119-1}.
\bibitem[{Russell and Schworm(2011)}]{russell.schworm.11}
\bibinfo{author}{Russell, R.R.}, \bibinfo{author}{Schworm, W.},
  \bibinfo{year}{2011}.
\newblock \bibinfo{title}{Properties of inefficiency indexes on $\langle$input,
  output$\rangle$ space}.
\newblock \bibinfo{journal}{J. Prod. Anal.} \bibinfo{volume}{36},
  \bibinfo{pages}{143--156}.
\newblock \DOIprefix\doi{10.1007/s11123-011-0209-3}.
\bibitem[{Russell and Schworm(2018)}]{russell.schworm.18}
\bibinfo{author}{Russell, R.R.}, \bibinfo{author}{Schworm, W.},
  \bibinfo{year}{2018}.
\newblock \bibinfo{title}{Technological inefficiency indexes: {A} binary
  taxonomy and a generic theorem}.
\newblock \bibinfo{journal}{J. Prod. Anal.} \bibinfo{volume}{49},
  \bibinfo{pages}{17--23}.
\newblock \DOIprefix\doi{10.1007/s11123-017-0518-2}.
\bibitem[{Sahoo et~al.(2014)Sahoo, Mehdiloozad and Tone}]{sahoo.al.14}
\bibinfo{author}{Sahoo, B.K.}, \bibinfo{author}{Mehdiloozad, M.},
  \bibinfo{author}{Tone, K.}, \bibinfo{year}{2014}.
\newblock \bibinfo{title}{Cost, revenue and profit efficiency measurement in
  {DEA}: a directional distance function approach}.
\newblock \bibinfo{journal}{European J. Oper. Res.} \bibinfo{volume}{237},
  \bibinfo{pages}{921--931}.
\newblock \DOIprefix\doi{10.1016/j.ejor.2014.02.017}.
\bibitem[{{\v{S}}ev\v{c}ovi\v{c} et~al.(2001){\v{S}}ev\v{c}ovi\v{c},
  Halick\'{a} and Brunovsk\'{y}}]{sevcovic.al.01}
\bibinfo{author}{{\v{S}}ev\v{c}ovi\v{c}, D.}, \bibinfo{author}{Halick\'{a},
  M.}, \bibinfo{author}{Brunovsk\'{y}, P.}, \bibinfo{year}{2001}.
\newblock \bibinfo{title}{D{EA} analysis for a large structured bank branch
  network}.
\newblock \bibinfo{journal}{Cent. Eur. J. Oper. Res.} \bibinfo{volume}{9},
  \bibinfo{pages}{329--342}.
\bibitem[{Shephard(1953)}]{shephard.53}
\bibinfo{author}{Shephard, R.W.}, \bibinfo{year}{1953}.
\newblock \bibinfo{title}{Cost and Production Functions}.
\newblock \bibinfo{publisher}{Princeton University Press, New Jersey}.
\bibitem[{Sueyoshi and Sekitani(2009)}]{sueyoshi.sekitani.09.ejor.196}
\bibinfo{author}{Sueyoshi, T.}, \bibinfo{author}{Sekitani, K.},
  \bibinfo{year}{2009}.
\newblock \bibinfo{title}{An occurrence of multiple projections in {DEA}-based
  measurement of technical efficiency: {T}heoretical comparison among {DEA}
  models from desirable properties}.
\newblock \bibinfo{journal}{European J. Oper. Res.} \bibinfo{volume}{196},
  \bibinfo{pages}{764--794}.
\newblock \DOIprefix\doi{10.1016/j.ejor.2008.01.045}.
\bibitem[{Tone(2001)}]{tone.01}
\bibinfo{author}{Tone, K.}, \bibinfo{year}{2001}.
\newblock \bibinfo{title}{A slacks-based measure of efficiency in data
  envelopment analysis}.
\newblock \bibinfo{journal}{European J. Oper. Res.} \bibinfo{volume}{130},
  \bibinfo{pages}{498--509}.
\newblock \DOIprefix\doi{10.1016/S0377-2217(99)00407-5}.
\bibitem[{Tone et~al.(2020)Tone, Chang and Wu}]{tone.al.20}
\bibinfo{author}{Tone, K.}, \bibinfo{author}{Chang, T.S.}, \bibinfo{author}{Wu,
  C.H.}, \bibinfo{year}{2020}.
\newblock \bibinfo{title}{Handling negative data in slacks-based measure data
  envelopment analysis models}.
\newblock \bibinfo{journal}{European J. Oper. Res.} \bibinfo{volume}{282},
  \bibinfo{pages}{926--935}.
\newblock \DOIprefix\doi{10.1016/j.ejor.2019.09.055}.

\end{thebibliography}
}
\end{document}